\documentclass[12pt]{article}
\usepackage{amssymb,amsfonts,latexsym,epic,eepic}
\title{{\bf Longo-Rehren subfactors arising from $\alpha$-induction}}
\author{{\sc Jens B\"ockenhauer} and {\sc David E. Evans}\\
School of Mathematics\\
University of Wales, Cardiff\\
PO Box 926, Senghennydd Road\\
Cardiff CF24 4YH, Wales, U.K.\\
e-mail: {\tt BockenhauerJM@cf.ac.uk, EvansDE@cf.ac.uk}  \\
\vphantom{X}\\
{\sc Yasuyuki Kawahigashi}\\
Department of Mathematical Sciences\\
University of Tokyo, Komaba, Tokyo, 153-8914, JAPAN\\
e-mail: {\tt yasuyuki@ms.u-tokyo.ac.jp}}
\date{February 18, 2000}

\begin{document}
\maketitle

\input amssym.def
\newsymbol\rtimes 226F

\def\Ad            {{\mathrm{Ad}}}
\def\Aut           {{\mathrm{Aut}}}
\def\bbC           {\mathbb{C}}
\def\bbM           {\mathbb{M}}
\def\bbN           {\mathbb{N}}
\def\bbNo          {\mathbb{N}_0}
\def\bbR           {\mathbb{R}}
\def\bbT           {\mathbb{T}}
\def\bbZ           {\mathbb{Z}}
\def\Z             {\mathbb{Z}}
\def\be            {\begin{equation}}
\def\bearl         {\begin{array}{l}}
\def\bearll        {\begin{array}{ll}}
\def\bearlll       {\begin{array}{lll}}
\def\bearrl        {\begin{array}{rl}}
\def\bea           {\begin{eqnarray}}
\def\beaa          {\begin{eqnarray*}}
\def\bfe           {{\bf1}}
\def\can           {\gamma}
\def\canr          {\theta}
\def\coioi         {{\co\iota^{\,-1}}}
\def\cA            {{\mathcal{A}}}
\def\cC            {{\mathcal{C}}}
\def\cCA           {\frak{C}}
\def\cD            {{\mathcal{D}}}
\def\cE            {{\mathcal{E}}}
\def\cF            {{\mathcal{F}}}
\def\cG            {{\mathcal{G}}}
\def\cH            {{\mathcal{H}}}
\def\cJ            {{\mathcal{J}}}
\def\cK            {{\mathcal{K}}}
\def\cL            {{\mathcal{L}}}
\def\cM            {{\mathcal{M}}}
\def\cN            {{\mathcal{N}}}
\def\cO            {{\mathcal{O}}}
\def\cP            {{\mathcal{P}}}
\def\cR            {{\mathcal{R}}}
\def\cS            {{\mathcal{S}}}
\def\cT            {{\mathcal{T}}}
\def\cV            {{\mathcal{V}}}
\def\cW            {{\mathcal{W}}}
\def\cX            {{\mathcal{X}}}
\def\cY            {{\mathcal{Y}}}
\def\cZ            {{\mathcal{Z}}}
\newcommand\co[1]  {\bar{{#1}}}
\def\diag          {{\mathrm{diag}}}
\def\dim           {{\mathrm{dim}}}
\newcommand\del[2] {\delta_{{#1},{#2}}}
\def\Dsys          {{\cD(\sys)}}
\def\E             {{\mathrm{e}}}
\def\ee            {\end{equation}}
\def\eear          {\end{array}}
\def\eea           {\end{eqnarray}}
\def\eeaa          {\end{eqnarray*}}
\def\End           {{\mathrm{End}}}
\def\eps           {\varepsilon}
\def\Eps           {{\mathcal{E}}}
\newcommand\erf[1] {Eq.\ (\ref{#1})}
\def\Exp           {{\mathrm{Exp}}}
\def\ext           {{\mathrm{ext}}}
\def\fh            {{\mathfrak{h}}}
\def\Gtwo          {{\mathrm{G}}_2}
\def\Hom           {{\mathrm{Hom}}}
\def\I             {{\mathrm{i}}}
\def\id            {{\mathrm{id}}}
\def\iotab         {{\co\iota}}
\def\Jz            {\mathcal{J}_z}
\def\lan           {\langle}
\def\lab           {{\co \lambda}}
\def\LG            {{\mathit{LG}}}
\def\LH            {{\mathit{LH}}}
\def\LIG           {{\mathit{L}}_I{\mathit{G}}}
\def\LIcG          {{\mathit{L}}_{\Ic}{\mathit{G}}}
\def\LIH           {{\mathit{L}}_I{\mathit{H}}}
\def\LIcSUn        {{\mathit{L}}_{I'}{\mathit{SU}}(n)}
\def\LISUn         {{\mathit{L}}_I{\mathit{SU}}(n)}
\def\LISUk         {{\mathit{L}}_I{\mathit{SU}}(k)}
\def\LISUnk        {{\mathit{L}}_I{\mathit{SU}}(nk)}
\def\LIcSUz        {{\mathit{L}}_{I'}{\mathit{SU}}(2)}
\def\LISUz         {{\mathit{L}}_I{\mathit{SU}}(2)}
\newcommand\ls[1]  {[\lambda_{{#1}}]}
\def\LSE           {L^2(S^1)}
\def\LSn           {L^2(S^1;\mathbb{C}^n)}
\def\LSUd          {{\mathit{LSU}}(3)}
\def\LSUk          {{\mathit{LSU}}(k)}
\def\LSUn          {{\mathit{LSU}}(n)}
\def\LSUnk         {{\mathit{LSU}}(nk)}
\def\LSUz          {{\mathit{LSU}}(2)}
\def\Mat           {{\mathrm{Mat}}}
\def\Mor           {{\mathrm{Mor}}}
\def\mub           {{\co{\mu}}}
\def\mult          {{\mathrm{mult}}}
\def\MXN           {{}_M {\cX}_N}
\def\MXM           {{}_M {\cX}_M}
\def\MXMa          {{}_M^{} {\cX}_M^\a}
\def\MXMo          {{}_M^{} {\cX}_M^0}
\def\MXMp          {{}_M^{} {\cX}_M^+}
\def\MXMm          {{}_M^{} {\cX}_M^-}
\def\MXMpm         {{}_M^{} {\cX}_M^\pm}
\def\NXN           {{}_N {\cX}_N}
\def\NXNd          {{}_N^{} {\cX}_N^{\mathrm{deg}}}
\def\NXM           {{}_N {\cX}_M}
\def\Nres          {\tilde{N}}
\def\NXM           {{}_N {\cX}_M}
\def\MYM           {{}_M {\cY}_M}
\def\MYMa          {{}_M^{} {\cY}_M^\a}
\def\MYMo          {{}_M^{} {\cY}_M^0}
\def\MYMp          {{}_M^{} {\cY}_M^+}
\def\MYMm          {{}_M^{} {\cY}_M^-}
\def\MYMpm         {{}_M^{} {\cY}_M^\pm}
\def\NYN           {{}_N {\cY}_N}
\def\NYNd          {{}_N^{} {\cY}_N^{\mathrm{deg}}}
\def\NYNper        {{}_N^{} {\cY}_N^{\mathrm{per}}}
\def\nub           {\overline{\nu}}
\def\op            {{\mathrm{opp}}}
\def\oto           {=0,1,2,\ldots}
\def\pio           {\pi_0}
\def\PSLZ          {{\mathit{PSL}}(2;\bbZ)}
\def\PSU           {{\mathit{PSU}}(1,1)}
\def\reso          {|_{\cH_0}}
\def\rmA           {{\mathrm{A}}}
\def\rmD           {{\mathrm{D}}}
\def\rmE           {{\mathrm{E}}}
\def\rmr           {{\mathrm{r}}}
\def\RXR           {{}_R {\cX}_R}
\def\Sect          {{\mathrm{Sect}}}
\def\SLC           {{\mathit{SL}}(2;\bbC)}
\def\SLnC          {{\mathit{SL}}(n;\bbC)}
\def\SLZ           {{\mathit{SL}}(2;\bbZ)}
\def\span          {{\mathrm{span}}}
\def\Ssys          {{\Sigma(\sys)}}
\def\SOf           {{\mathit{SO}}(5)}
\def\SON           {{\mathit{SO}}(N)}
\def\SUd           {{\mathit{SU}}(3)}
\def\SUk           {{\mathit{SU}}(k)}
\def\SUm           {{\mathit{SU}}(m)}
\def\SUn           {{\mathit{SU}}(n)}
\def\SUnk          {{\mathit{SU}}(nk)}
\def\SUz           {{\mathit{SU}}(2)}
\def\SUzk          {{\mathit{SU}}(2k)}
\def\SUf           {{\mathit{SU}}(4)}
\def\sys           {{\Delta}}
\def\tr            {{\mathrm{tr}}}
\def\Tr            {{\mathrm{Tr}}}
\def\Un            {\mathit{U}(n)}


\def\qed{{\unskip\nobreak\hfil\penalty50
\hskip2em\hbox{}\nobreak\hfil  $\Box$      
\parfillskip=0pt \finalhyphendemerits=0\par}\medskip}
\def\proof{\trivlist \item[\hskip \labelsep{\it Proof.\ }]}
\def\endproof{\null\hfill\qed\endtrivlist}

\def\equi{\sim}
\def\isom{\cong}
\def\ti{\tilde}
\def\lan{\langle}
\def\ran{\rangle}
\def\a{\alpha}
\def\de{\delta}
\def\ga{\gamma}
\def\Ga{\Gamma}
\def\la{\lambda}
\def\La{\Lambda}
\def\th{\theta}
\def\om{\omega}
\def\Om{\Omega}
\def\Ups{\Upsilon}
\def\si{\sigma}

\newcommand\labl[1]{\label{#1}}


\def\thinlines{\allinethickness{0.3pt}}
\def\thicklines{\allinethickness{1.0pt}}
\def\Thicklines{\allinethickness{2.0pt}}


\newtheorem{theorem}{Theorem}[section]
\newtheorem{lemma}[theorem]{Lemma}
\newtheorem{conjecture}[theorem]{Conjecture}
\newtheorem{corollary}[theorem]{Corollary}
\newtheorem{definition}[theorem]{Definition}
\newtheorem{assumption}[theorem]{Assumption}
\newtheorem{proposition}[theorem]{Proposition}
\newtheorem{remark}[theorem]{Remark}
\newtheorem{example}[theorem]{Example}

\begin{abstract}
We study (dual) Longo-Rehren subfactors
$M\otimes M^\op \subset R$ arising from various systems
of endomorphisms of $M$ obtained from $\a$-induction for
some braided subfactor $N\subset M$.
Our analysis provides useful tools to determine the systems
of $R$-$R$ morphisms associated with such Longo-Rehren
subfactors, which constitute the ``quantum double'' systems
in an appropriate sense.
The key to our analysis is that $\a$-induction
produces half-braidings in the sense of Izumi,
so that his general theory can be applied.
Nevertheless, $\a$-induced systems are in general
not braided, and thus our results allow to compute the
quantum doubles of (certain) systems without braiding.
We illustrate our general results by several examples,
including the computation of the quantum double systems
for the asymptotic inclusion of the E$_8$ subfactor as well
as its three analogues arising from conformal inclusions
of $SU(3)_k$.
\end{abstract}

\newpage

\section{Introduction}

There are various constructions analogous to the quantum double
construction of Drinfel$'$d \cite{D} in subfactor theory. 
The first of such constructions is Ocneanu's
asymptotic inclusion (see e.g.\ \cite{EK2}) which produces
$M\vee (M'\cap M_\infty)\subset M_\infty$ from a
given hyperfinite II$_1$ subfactor $N\subset M$ with finite
index and finite depth.  That is, if we compare the system
of $M$-$M$ bimodules (or $N$-$N$ bimodules) arising from
$N\subset M$ and that of $M_\infty$-$M_\infty$ bimodules
arising from $M\vee (M'\cap M_\infty)\subset M_\infty$,
then the latter can be regarded as a ``quantum double'' of the
former due to its categorical structure.
This viewpoint was noticed by Ocneanu in connection
to topological quantum field theory of three dimensions, and
the categorical meaning of the construction has been 
recently clarified by M\"uger \cite{Mu}.
(A general reference for the asymptotic inclusions and
topological quantum field theories is \cite[Chapter 12]{EK2}.
We actually need a connectedness assumption of a certain
graph for the above interpretation of ``quantum double.''
See \cite[Theorem 12.29]{EK2} for a precise statement.)
Popa's notion of a symmetric enveloping algebra in \cite{P}
also gives a construction of a new subfactor from a given one,
and if  the initial subfactor $N\subset M$ is hyperfinite,
of type II$_1$, of finite index, and of finite depth,
then this construction gives a subfactor isomorphic to
the asymptotic inclusion.

Later, Longo and Rehren introduced in \cite{LR} another
construction of a subfactor from a given system $\sys$ of
endomorphisms, which is now called the Longo-Rehren subfactor.
Masuda \cite{M} has proved that the asymptotic inclusion and
the Longo-Rehren subfactor are essentially
the same constructions, though the constructions arise from
very different viewpoints and appear rather unrelated
at first sight.
Izumi \cite{I1} has developed a general theory on the
structure of sectors associated with Longo-Rehren subfactors.
He introduced a notion of half-braiding and showed
that the structure of the quantum double system $\Dsys$
is closely related to half-braidings. Namely,
any morphism in $\Dsys$ is given by certain extensions
of morphisms defined by means of a half-braiding.
This extension will be called $\eta$-extension in this paper.
Moreover, he presents various interesting applications
in \cite{I2} with calculations involving Ocneanu's tube
algebra handled in the setting of Longo-Rehren subfactors.

Longo and Rehren also introduced an extension formula for
endomorphisms of a smaller net to a larger net for nets of
subfactors in the same paper \cite{LR}. Xu \cite{X1,X2}
obtained various interesting results by using essentially
the same construction in connection to conformal inclusions.
Two of us \cite{BE1,BE3} systematically analyzed the
extension formula of Longo and Rehren for nets of subfactors.
It was named $\a$-induction in \cite{BE1,BE2,BE3} in order
to emphasize structural similarities with the Mackey machinery
of induction and restriction of group representations and
to distinguish it from the different sector induction,
nevertheless.
We have further studied $\a$-induction in the very general
setting of braided subfactors in \cite{BEK1,BEK2}.
We identified it with Ocneanu's graphical construction of
chiral generators and obtained several results by making
use of his graphical methods of double triangle algebras.
Izumi's work \cite{I1} shows that the study of Longo-Rehren
subfactors using a half-braiding is somewhat similar to the
study of $\a$-induction. Moreover, $\a$-induction produces
interesting systems of endomorphisms which come with
various half-braidings, as we will demonstrate in this paper.
So it is quite natural to study their quantum doubles by means
of associated Longo-Rehren subfactors and applying Izumi's
general theory, and this is what we propose in this paper.

To be more specific, we start with a subfactor $N\subset M$
with a finite braided system of $N$-$N$ morphisms allowing
us to apply $\a$-induction.
Then the two chiral $\a$-inductions arising from the
braiding produce chiral systems of $M$-$M$ morphisms,
and together they generate the full induced system.
We define a system $\sys$ to be (subsystems of) either the
chiral or the full induced system and study their associated
Longo Rehren inclusions $M\otimes M^\op \subset R(\sys)$.
We construct half-braidings with respect to such systems $\sys$
for certain classes of endomorphisms, and this enables us to
apply Izumi's theory for analyzing the structure of
the quantum double system $\Dsys$ which can be
given by $\eta$-extensions.
The important point is that $\a$-induced systems are
not braided in general. (They can even be
non-commutative \cite{X1,X2,BE2,BE3} and general criteria
for non-commutativity were given in \cite{BEK1,BEK2}.)
Thus our analysis is aimed at going beyond the
computation of quantum doubles of braided systems
which has been carried out in \cite{O6,EK1},
and avoiding at the same time complex constructions
like the tube algebra used in \cite{I2}.

In fact, the rich structure of $\a$-induced systems
allows to derive fairly concrete results concerning
the structure of the quantum doubles.
Namely, we derive concrete formulae for the
(dimensions of the) intertwiner spaces between
various $\eta$-extensions.
It is crucial that we allow the braiding on the
$N$-$N$ morphisms to be degenerate.
However, the situation simplifies considerably
whenever this braiding is non-degenerate.
For example, in this case the quantum double of the
full induced system is given as the direct product
of the original $N$-$N$ system with itself.
As a corollary we obtain a new proof of Rehren's recent
theorem on ``generalized Longo-Rehren subfactors'' in a
typical case arising from $\a$-induction.
Similarly, the quantum double of the chiral system
is given by the direct product of the original
$N$-$N$ system with the ambichiral system in the
non-degenerate case.
However, in the general, degenerate case the situation
is more involved. More precisely, the subsystem of
degenerate morphisms arranges the direct product of the
$N$-$N$ system with the ambichiral system into orbits whose
elements have to be identified whereas fixed points split,
so that the quantum double is now some kind of orbifold of
the one we would have obtained in the non-degenerate case.

An orbifold phenomenon has been encountered
earlier in computations of dual principal graphs and
bimodule systems of asymptotic inclusions of $\SUn_k$
subfactors which correspond to degenerately braided
systems \cite{O6,EK1}. Our results show that
the same phenomenon shows up for quantum doubles of
(in general not braided) systems arising from $\a$-induction
and having its origin in degeneracies of the braiding
of the original $N$-$N$ system.
We illustrate this by computing the quantum doubles
of several examples arising from conformal inclusions
of $\SUz$ and $\SUd$.
They correspond to the asymptotic inclusions of subfactors
with principal graph E$_6$ and E$_8$ as well their three
analogues from $\SUd$ conformal inclusions.

This paper is organized as follows. In Section \ref{prelim} we
recall basic facts on $\a$-induc\-tion, state our main assumption
and review the results of \cite{I1} we use in the sequel.
In Section \ref{sectMXMa} we
consider the quantum doubles of full induced systems.
We introduce half-braidings for the induced morphisms
$\a^\pm_\la$ and obtain formulae for the intertwiner spaces
of their $\eta$-extensions and finally consider the
non-degenerate case.
In Section \ref{sectMXMpm} we propose the same analysis
for the quantum doubles of chiral systems.
Finally we treat examples arising from conformal
inclusions in Section \ref{col0}.

\section{Preliminaries}
\labl{prelim}
\subsection{Braided systems of morphisms and $\a$-induction}

Let $A$ and $B$ be type III von Neumann factors.
A unital $\ast$-homomorphism $\rho:A\rightarrow B$
is called a $B$-$A$ morphism. The positive number
$d_\rho=[B:\rho(A)]^{1/2}$ is called the statistical
dimension of $\rho$; here $[B:\rho(A)]$ is the
Jones index \cite{J} of the subfactor
$\rho(A)\subset B$. If $\rho$ and $\sigma$ are $B$-$A$
morphisms with finite statistical dimensions, then
the vector space of intertwiners
\[ \Hom(\rho,\sigma)=\{ t\in B: t\rho(a)=\sigma(a)t \,,
\,\, a\in A \}  \]
is finite-dimensional, and we denote its dimension by
$\lan\rho,\sigma\ran$.
An $A$-$B$ morphism $\co\rho$ is a
conjugate morphism if there are isometries
$r_\rho\in\Hom(\id_A,\co\rho\rho)$ and
${\co r}_\rho\in\Hom(\id_B,\rho\co\rho)$ such that
$\rho(r_\rho)^* {\co r}_\rho=d_\rho^{-1}\bfe_B$ and
$\co\rho({\co r}_\rho)^* r_\rho=d_\rho^{-1}\bfe_A$.
The map $\phi_\rho:B\rightarrow A$,
$b\mapsto r_\rho^* \co\rho(b)r_\rho$, is called the
(unique) standard left inverse and satisfies
\be
\phi_\rho(\rho(a)b\rho(a'))=a\phi_\rho(b)a' \,,
\quad a,a'\in A\,,\quad b\in B \,.
\label{phiaba}
\ee

We work with the setting of \cite{BEK1}, i.e.\ we are
working with a type III subfactor and
finite system $\NXN\subset\End(N)$ of
braided morphisms which is compatible with the inclusion
$N\subset M$. Then the inclusion is in particular forced to have
finite Jones index and also finite depth (see e.g.\ \cite{EK2}).
More precisely, we make the following

\begin{assumption}
{\rm Let $N\subset M$ be a type III subfactor
together with a finite system of endomorphisms
$\NXN\subset\End(N)$ in the sense of \cite[Def.\ 2.1]{BEK1}
which is braided in the sense of \cite[Def.\ 2.2]{BEK1}.
For a given subsystem $\NYN\subset\NXN$ we assume that
$\canr=\co\iota\iota\in\Sigma(\NYN)$ for the
injection $M$-$N$ morphism $\iota:N\hookrightarrow M$ and a
conjugate $N$-$M$ morphism $\co\iota$.}
\labl{assump}
\end{assumption}
Here $\Sigma(\NYN)$ denotes the set of finite sums of
morphisms in $\NYN$, and we will use a similar notation
for other systems.

With the braiding $\eps$ on $\NXN$ and its
extension to $\Sigma(\NXN)$ as in \cite{BEK1}, one can
define the $\a$-induced morphisms $\a^\pm_\la\in\End(M)$
for $\la\in\Sigma(\NXN)$ by the Longo-Rehren formula \cite{LR},
namely by putting
\[ \a_\la^\pm = \co\iota^{\,-1} \circ \Ad
(\eps^\pm(\lambda,\canr)) \circ \lambda \circ \co\iota \,, \]
where $\co\iota$ denotes a conjugate morphism of the
injection map $\iota:N\hookrightarrow M$.
Then $\a^+_\la$ and $\a^-_\la$ extend $\la$, i.e.\
$\a^\pm_\la\circ\iota=\iota\circ\la$, which in turn implies
$d_{\a_\la^\pm}=d_\la$ by the multiplicativity of
the minimal index \cite{L3}. 
Let $\can=\iota\co\iota$ denote Longo's canonical
endomorphism from $M$ into $N$. Then there is an isometry
$v\in\Hom(\id,\can)$ such that any $m\in M$ is uniquely
decomposed as $m=nv$ with $n\in N$.
Thus the action of the extensions
$\a^\pm_\la$ are uniquely characterized by the relation
$\a^\pm_\la(v)=\eps^\pm(\la,\canr)^* v$ which can be
derived from the braiding fusion equations
(BFE's, see e.g.\ \cite[Eq.\ (5)]{BEK1}).
Moreover, we have $\a_{\la\mu}^\pm=\a_\la^\pm \a_\mu^\pm$
if also $\mu\in\Sigma(\NXN)$, and clearly
$\a_{{\mathrm{id}}_N}^\pm={{\mathrm{id}}}_M$.
In general one has
\be
\Hom(\la,\mu) \subset \Hom(\a^\pm_\la,\a^\pm_\mu)
\subset \Hom(\iota\la,\iota\mu) \,, \qquad \la,\mu\in\Sigma(\NXN) \,.
\label{ahomincl}
\ee
The first inclusion is a consequence of the BFE's.
Namely, $t\in\Hom(\la,\mu)$ obeys
$t\eps^\pm(\canr,\la)=\eps^\pm(\canr,\mu)\canr(t)$,
and thus
\[ t \a_\la^\pm(v) = t \eps^\pm(\la,\canr)^* v =
\eps^\pm(\mu,\canr)^* \canr(t) v
=\eps^\pm(\mu,\canr)^* vt = \a_\mu^\pm(v) t \,. \]
The second follows from the extension property
of $\a$-induction. Hence $\a_{\co\la}^\pm$
is a conjugate for $\a_\la^\pm$ as there are
$r_\la\in\Hom(\id,\co\la\la)\subset
\Hom(\id,\a_{\co\la}^\pm\a_\la^\pm)$
and ${\co r}_\la\in\Hom(\id,\la\co\la)\subset
\Hom(\id,\a_\la^\pm\a_{\co\la}^\pm)$ such
that
$\la(r_\la)^*{\co r}_\la=\co\la({\co r}_\la)^*r_\la
= d_\la^{-1}\bfe$. We also have some kind of
naturality equations for $\a$-induced morphisms,
\begin{equation}
x \eps^\pm(\rho,\la)=\eps^\pm(\rho,\mu)\a^\pm_\rho(x)
\label{anat}
\ee
whenever $x\in\Hom(\iota\la,\iota\mu)$,
$\rho\in\Sigma(\NXN)$.

Recall that the statistics phase of $\om_\la$ for
$\la\in\NXN$ is given as
\[ d_\la \phi_\la(\eps^+(\la,\la))=\om_\la \bfe \,. \]
The monodromy matrix $Y$ is defined by
\be
Y_{\la,\mu} = \sum_{\rho\in\NXN}
\frac{\om_\la \om_\mu}{\om_\rho} N_{\la,\mu}^\rho d_\rho \,,
\qquad \la,\mu\in\NXN \,,
\label{Ymat}
\ee
with $N_{\la,\mu}^\rho=\lan\rho,\la\mu\ran$ denoting the
fusion coefficients. Then one checks that $Y$ is symmetric,
that $Y_{\co\la,\mu}=Y_{\la,\mu}^*$ as well as
$Y_{\la,0}=d_\la$ \cite{R1,FG,FRS2}. (As usual, the
label ``$0$'' refers to the identity morphism $\id\in\NXN$.)
Now let $\Omega$ be the diagonal matrix with entries
$\Omega_{\la,\mu}=\om_\la \del\la\mu$.
Putting
\be
Z_{\la,\mu}=\lan\a^+_\la,\a^-_\mu\ran\,,\qquad \la,\mu\in\NXN\,,
\label{couplmat}
\ee
defines a matrix subject to the constraints
\[ Z_{\la,\mu}=0,1,2,\ldots \,, \qquad
\mbox{and} \qquad Z_{0,0}=1 \,, \]
and commuting with $Y$ and $\Omega$ \cite{BEK1}.
The Y- and $\Omega$-matrices obey
$\Omega Y \Omega Y \Omega = z Y$
where $z=\sum_\la d_\la^2 \om_\la$ \cite{R1,FG,FRS2},
and this actually holds even if the braiding is
degenerate (see \cite[Sect.\ 2]{BEK1}).
If $z\neq 0$ we put $c=4\arg(z)/\pi$, which is defined
modulo 8, and call it the ``central charge''.
Moreover, S- and T-matrices are then defined by
\[ S = |z|^{-1} Y \,, \qquad T = \E^{-\I\pi c/12} \Omega \]
and hence fulfill $TSTST=S$. 
One has $|z|^2=[[\NXN]]$ with the global index
$[[\NXN]]=\sum_\la d_\la^2$ and $S$ is unitary, so that
$S$ and $T$ are indeed the standard generators in a
unitary representation of the modular group $\SLZ$,
if and only if the braiding is non-degenerate \cite{R1}.
Consequently, $Z$ gives a modular invariant in this case.

Let $\MXM\subset\End(M)$ denote a system of endomorphisms
consisting of a choice of representative endomorphisms of
each irreducible subsector of sectors of the form
$[\iota\la\co\iota]$, $\la\in\NXN$.
We choose $\id\in\End(M)$
representing the trivial sector in $\MXM$.
Then we define similarly the chiral systems
$\MXMpm$ and the $\a$-system $\MXMa$ to be the subsystems
of endomorphisms $\beta\in\MXM$
such that $[\beta]$ is a subsector of $[\a^\pm_\la]$ and of
of $[\a_\la^+\a_\mu^-]$, respectively,
for some $\la,\mu\in\NXN$.
(Note that any subsector of $[\a_\la^+\a_\mu^-]$ is
automatically a subsector of  $[\iota\nu\co\iota]$
for some $\nu\in\NXN$.)
The ambichiral system is defined
as the intersection $\MXMo=\MXMp\cap\MXMm$, so that
$\MXMo \subset \MXMpm \subset \MXMa \subset \MXM$.
Thus their ``global indices'', i.e.\ the sums over the squares
of the statistical dimensions of their morphisms, fulfill
$1 \le [[\MXMo]] \le [[\MXMpm]] \le [[\MXMa]] \le [[\MXM]]=[[\NXN]]$.
(Throughout this paper we denote the global index of a
system by use of double rectangular brackets.)

Let us now consider the subsystem $\NYN$ appearing in
Assumption \ref{assump}. If the inclusion $\NYN\subset\NXN$
is proper, then we may play the same game considering
$\a$-induction for exclusively $\la\in\NYN$.
This way we will obtain $\a$-induced
systems which are contained in the $\a$-induced systems
associated to $\NXN$, i.e.\ we have the following
scheme of inclusions:
\[ \begin{array}{ccccccc}
\MXMo & \subset & \MXMpm & \subset & \MXMa & \subset & \MXM \\
 \cup &         &  \cup  &         &  \cup &         & \cup \\
\MYMo & \subset & \MYMpm & \subset & \MYMa & \subset & \MYM
\end{array} \]
We will use these systems for the construction of
Longo-Rehren subfactors and for the analysis of sectors
associated to them.
We are particularly interested in examples where
(at least) the braiding on the subsystem $\NYN$ may be
degenerate. Let $\NYNd$ denote the system of degenerate
morphisms, i.e.\
\[ \NYNd = \{ \nu\in\NYN \mid \eps^+(\nu,\rho)=\eps^-(\nu,\rho)
\quad \mbox{for all} \quad \rho\in\NYN  \} \,. \]
Clearly, the braiding on $\NYN$ is non-degenerate
(in the sense of \cite{R1} or \cite[Def.\ 2.3]{BEK1})
if and only if $\NYNd=\{\id\}$.
Note that since $\canr$ decomposes by Assumption \ref{assump}
only into morphisms of $\NYN$ and since
$\a^\pm_\la(v)=\eps^\pm(\la,\canr)^*v$ for any $\la\in\NXN$
we find $\a^+_\rho=\a^-_\rho$ whenever $\rho\in\NYNd$.
Finally we introduce
\[ \NYNper = \{ \la\in\NXN \mid \eps^+(\la,\rho)=\eps^-(\la,\rho)
\quad \mbox{for all} \quad \rho\in\NYN \} \]
and call it the {\sl relative permutant} of $\NYN$ in $\NXN$.
Clearly, $\a^+_\la=\a^-_\la$ whenever $\la\in\NYNper$.

\subsection{Longo-Rehren subfactors,
half-braidings and $\eta$-ex\-ten\-sions}

Let $M$ be a type III factor with a finite system
$\sys\subset\End(M)$ of endomorphisms. Let
$M^\op$ denote the opposite algebra of $M$ and consider
$M\otimes M^\op$.
By constructing a ``Q-system'', Longo and Rehren
showed in \cite[Prop.\ 4.10]{LR} that there
is a (type III) subfactor $B \subset M\otimes M^\op$
with canonical endomorphism $\Theta\in\End(M\otimes M^\op)$
decomposing as a sector as
\[ [\Theta] =\bigoplus_{\beta\in\sys}
\,\,\,[\beta\otimes\beta^\op] \,. \]
Here $\beta^\op=j\circ\beta\circ j^{-1}$ where
$j:M\rightarrow M^\op$ is the anti-linear isomorphism.
The subfactor $B \subset M\otimes M^\op$ is now
called the Longo-Rehren subfactor.
For reasons of convenience, we consider in this paper
the dual subfactor $M\otimes M^\op\subset R$ and call it
the Longo-Rehren subfactor as well.
(This convention is compatible with \cite{KLM}.)
That is, $B \subset M\otimes M^\op\subset R$ is a Jones
extension and $\Theta$ is then the dual canonical
endomorphism of $M\otimes M^\op\subset R$.

The following is a slight variation of Izumi's definition
\cite[Def.\ 4.2]{I1} of a half-braiding.

\begin{definition}{\rm 
Let $\Phi$ be a system of morphisms in $\End(M)$ and
$\sys\subset\Phi$ a subsystem.
For $\sigma\in\Sigma(\Phi)$ we call a family of
unitary operators
$\Eps_\si=\{\Eps_\si(\beta)\}_{\beta\in \sys}$ 
a {\sl half-braiding} of $\sigma$ with respect to
$\sys$ if it satisfies the following two conditions:
\begin{enumerate}
\item $\Eps_\si(\beta)\in \Hom(\si\beta, \beta\si)$
for all $\beta\in\sys$. 
\item  Whenever $\beta_1,\beta_2,\beta_3\in\sys$ then 
\[ X\Eps_\si(\beta_3)=\beta_1(\Eps_\si(\beta_2))
\Eps_\si(\beta_1)\si(X) \]
holds for every $X\in \Hom(\beta_3,\beta_1\beta_2)$,
\end{enumerate}
Two pairs $(\si,\Eps_\si)$, $(\si',\Eps'_{\si'})$
of morphisms $\si,\si'\in\Sigma(\Phi)$ with respective
half-braidings $\Eps_\si,\Eps'_{\si'}$ are said to be
{\sl equivalent} if there is unitary
$u\in \Hom(\sigma',\sigma)$ such that
\[\Eps_\si(\beta)=\beta(u)\Eps'_\si(\beta)u^*\]
for all $\beta\in\sys$.}
\labl{half}
\end{definition}

Note that our definition of equivalence is slightly more general 
than the one in \cite[Def.\ 4.2]{I1} because we choose the
$\sigma$'s from a generically larger set $\Phi\supset\sys$.
We then define an extension $\eta(\si,\Eps_\si)$ of the endomorphism
$\si\otimes\id$ of $M\otimes M^\op$ to $R$ as in the following
definition, which is just the dual version of
Izumi's definition of $(\widetilde{\si\!,\!\Eps_\si})$
in \cite[Def.\ 4.4]{I1}.
This extension is somewhat similar to $\a$-induction.
Izumi's important observation is that we need only ``half''
the properties of a usual braiding for this extension.
We need some preparation.
Let $W_\beta\in\Hom(\beta\otimes\beta^\op,\Theta)$,
$\beta\in\sys$, be isometries so that
$W_\beta^* W_{\beta'}=\delta_{\beta,\beta'}\bfe$
and $\sum_{\beta\in\sys}W_\beta W_\beta^*=\bfe$.
(Note that for a Longo-Rehren subfactor with given
$\Theta$ each $W_\beta$ is unique up to a phase.)
Let $\iota_{\mathrm{LR}}: M\otimes M^\op \hookrightarrow R$
denote the inclusion homomorphism so that the dual
canonical endomorphism is given by
$\Theta=\bar{\iota}_{\mathrm{LR}}\iota_{\mathrm{LR}}$,
and then
$\Gamma=\iota_{\mathrm{LR}}\bar{\iota}_{\mathrm{LR}}$
is a canonical endomorphism.
Then there is \cite{L2} an isometry
$V\in\Hom(\id,\Gamma)$ such that
$W_\id^*V=[R: M\otimes M^\op]^{-1/2}\bfe$, and note
that $[R: M\otimes M^\op]=\sum_{\beta\in\sys}d_\beta^2$.
Moreover, for each $X\in R$ there is a unique
$a\in M\otimes M^\op$ such that $X=aV$.

\begin{definition}{\rm
For $\si\in\Sigma(\Phi)$ with a half-braiding
$\Eps_\si=\{\Eps_\si(\beta)\}_{\beta\in\sys}$,
we define an extension $\eta(\si,\Eps_\si)\in\End(R)$
by putting
\be
\begin{array}{rl}
\eta(\si,\Eps_\si)(a)&=(\si\otimes\id)(a),
\quad a\in M\otimes M^\op,\\[.4em]
\eta(\si,\Eps_\si)(V)&=U(\si,\Eps_\si)^* V,
\end{array}
\ee
where the unitary $U(\si,\Eps_\si)$ is defined as
\be
U(\si, \Eps_\si)=\sum_{\beta\in\sys}
W_\beta(\Eps_\si(\beta)\otimes 1)
(\si\otimes \id^\op)(W_\beta^*) \,.
\ee
Using
\be
U^\op(\si, \Eps_\si)=\sum_\beta W_\beta(1\otimes
j(\Eps_\si(\beta)))(\id\otimes\si^\op)(W_\beta^*)
\ee
we similarly define an extension
$\eta^\op(\si,\Eps_\si)\in\End(R)$ of
$\id\otimes\si^\op$.}
\labl{eta}
\end{definition}

Let $\Dsys$ be the system of irreducible endomorphisms of $R$
arising from a choice of representative morphisms of
irreducible subsectors of
$\iota_{\mathrm{LR}}\circ\beta'\otimes
\beta^\op\circ\bar{\iota}_{\mathrm{LR}}$
for $\beta,\beta'\in\sys$.
Following \cite[Def.\ 4.4]{I1}, we call $\Dsys$ the
{\sl quantum double system} of $\sys$. 
(Note that Izumi's notation $D(\sys)$ for the
quantum double includes reducible morphisms and
thus corresponds to $\Sigma(\Dsys)$.
Also note that the system $\Dsys$ may be strictly
larger than that arising from the Longo-Rehren subfactor
$M\otimes M^\op \subset R$ in the usual sense,
i.e.\ arising from the decomposition of all powers of $\Gamma$.
See Remark after \cite[Thm.\ 4.6]{I1}.)
Izumi has proved in \cite[Lemma 4.5, Thm.\ 4.6]{I1}
that $\eta(\si,\Eps)$ gives an endomorphism in $\Sigma(\Dsys)$
if we consider $\si\in\Sigma(\sys)$ only, and then any
endomorphism in $\Sigma(\Dsys)$ arises in this way.
Note that this will no longer be true if we consider
generic $\si\in\Sigma(\Phi)$.

The following is nothing but Izumi's
\cite[Thm.\ 4.6 (ii)]{I1}. We only provide a proof
for the reader's convenience and in order to
demonstrate that the arguments are the same
though we work in a picture dual to Izumi's and
extend $\si\in\Sigma(\Phi)\supset\Sigma(\sys)$.

\begin{theorem}
Let $\si,\si'\in\Sigma(\Phi)$ with half-braidings
$\Eps_\si=\{\Eps_\si(\beta)\}_{\beta\in\sys}$,
$\Eps_{\si'}'=\{\Eps_{\si'}'(\beta)\}_{\beta\in\sys}$.
Then we have
\be
\begin{array}{l}
\Hom(\eta(\si,\Eps_\si),\eta(\si',\Eps_{\si'}'))=\\[.4em]
\qquad\quad=\{X\otimes\bfe\mid X\in \Hom(\si,\si')\,,
\,\,\,\Eps'_{\si'}(\beta)X= \beta(X)\Eps_\si(\beta)
{\rm \ for\ all\ }\beta\in \sys\}\,.
\end{array}
\ee
In particular, $\eta(\si,\Eps_\si)$ and
$\eta(\si',\Eps_{\si'}')$ are unitarily equivalent
as morphisms of $R$ if and only if pairs
$(\si,\Eps_\si)$ and $(\si',\Eps_{\si'}')$ are equivalent
in the sense of Definition \ref{half}.
\labl{irred}
\end{theorem}

\begin{proof}
Let $T\in\Hom(\eta(\si,\Eps_\si),\eta(\si',\Eps_{\si'}'))$.
Then it is decomposed as $T=aV$ with
$a\in\Hom(\Theta\circ(\si\otimes\id^\op),\si'\otimes\id^\op)$.
Consequently
$aW_\beta\in\Hom(\beta\si\otimes\beta^\op),\si'\otimes\id^\op)$
can be non-zero only for $\beta=\id$.
Hence $a=bW_\id^*$ with
\[ b=aW_\id\in\Hom(\si\otimes\id^\op,\si'\otimes\id^\op)
=\{X\otimes\bfe\mid X\in\Hom(\si,\si')\}\,.\]
Since $W_\id^*V$ is a (non-zero) scalar we have found
$T\in\{X\otimes\bfe\mid X\in\Hom(\si,\si')\}$.
For such a $T=X\otimes\bfe$ the condition
$TU(\si,\Eps_\si)^*V=U(\si',\Eps_{\si'}')^*VT$
is equivalent to
$\Theta(X\otimes\bfe)U(\si,\Eps_\si)
=U(\si',\Eps_{\si'}')X\otimes\bfe$.
Sandwiching with $W_\beta^*$ and
$(\sigma\otimes\id^\op)(W_\beta)$ gives the
desired intertwining relations for all $\beta\in\sys$.
Conversely, any $T=X\otimes\bfe$ with $X\in\Hom(\si,\si)$
satisfying these relations intertwines
$\eta(\si,\Eps_\si)$ and $\eta(\si',\Eps_{\si'}')$.
\end{proof}

Since $\eta(\si,\Eps_\si)$
is an extension of $\sigma\otimes\id^\op$ and since
$[R:M\otimes M^\op]<\infty$ we also find
that its statistical dimension is $d_\si$,
i.e.\ $\eta$ preserves statistical dimensions.
We have even more than that. Namely, for pairs
$(\si,\Eps_\si)$ as above, we have natural notions of
addition and multiplication extending those of the
endomorphisms $\si$. Let $\si_i\in\Sigma(\Phi)$
with half-braidings $\Eps^i_{\si_i}$,
$i=1,2,...,n$. Let $\{t_i\}_{i=1}^n$ be a set
of isometries in $M$ satisfying the Cuntz relations
and let $\si\in\Sigma(\Phi)$ be given by
$\si(m)=\sum_i t_i \si_i(m) t_i^*$ for all $m\in M$.
It is routine to show that putting
\[ \Eps_\si(\beta) = \sum_{i=1}^n \beta(t_i)
\Eps^i_{\si_i}(\beta)t_i^* \,, \qquad \beta\in\sys \,,\]
defines a half-braiding for $\sigma$.
Similarly, putting
\[ \Eps_{\si'}'(\beta)= \Eps^1_{\si_1}(\beta)
\si_1(\Eps^2_{\si_2}(\beta)) \,, \qquad \beta\in\sys \,,\]
defines a half-braiding
$\{\Eps_{\si'}'(\beta)\}_{\beta\in\sys}$ of
products $\si'=\si_1\si_2$, as used \cite{I1}.
It is straightforward to show that we have exact
multiplicativity for the $\eta$-extensions,
\[ \eta(\si',\Eps_{\si'})=\eta(\si_1,\Eps^1_{\si_1})
\eta(\si_2,\Eps^2_{\si_2}) \,, \]
with this product half-braiding.
Finally, conjugates were defined in \cite[Thm.\ 4.6 (iv)]{I1}
as follows. For a pair $(\si,\Eps_\si)$, operators
\[ \bar\Eps_{\bar\si}(\beta)=
d_\si R^*_\si \bar\si(\Eps_\si(\beta)^*
\beta(\bar R_\si)) \,, \qquad \beta\in\sys \,, \]
where $R_\si\in\Hom(\id,\bar\si \si)$,
$\bar R_\si\in\Hom(\id,\si \bar\si)$ are isometries with
$\bar R^*_\si \si(R_\si)= 
R^*_\si \bar\si(\bar R_\si)=d_{\si}^{-1}$,
give a half-braiding for the conjugate morphism $\bar\si$.
The half-braiding $\{\bar\Eps_{\bar\si}(\beta)\}_{\beta\in\sys}$
depends on the choices of $R_\si,\bar R_\si$ in general, however,
its equivalence class does not \cite{I1}.
Then Izumi's results give the following

\begin{proposition}
The extension map
$\eta:(\si,\Eps_\si)\rightarrow\eta(\si,\Eps_\si)$,
regarded as a map from equivalence classes of pairs
to sectors of $R$, preserves the operations of addition,
multiplication, and conjugates.
\labl{homom}
\end{proposition}

\begin{proof}
The preservation of addition and the multiplication
is a straight-forward corollary of Theorem \ref{irred}.
The statement for the conjugates is derived in the same
way as \cite[Thm.\ 4.6 (iv)]{I1}.
\end{proof}

Next, \cite[Prop.\ 6.4]{I1} gives the following

\begin{proposition}
For a pair $\si\in\Sigma(\Phi)$ with a half-braiding
$\{\Eps_\si(\beta)\}_{\beta\in\sys}$, the extensions
$\eta(\si,\Eps_\si)$ and
$\eta^\op(\bar\si,\bar\Eps_{\bar\si})$
are unitarily equivalent.
\labl{conj}
\end{proposition}

Finally, \cite[Thm.\ 4.1]{I1} and the remark at the end
of \cite[Sect.\ 4]{I1} give the following

\begin{proposition}
Let $\cG$ be the bipartite graph with odd vertices
labelled by $\sys$ and even vertices labelled $\Dsys$,
and the number of edges between a vertex labelled by
$\beta\in\sys$ and a vertex labelled by $\Omega\in\Dsys$
such that $[\Omega]=[\eta(\si,\Eps_\si)]$ for
$\si\in\Sigma(\sys)$ with some half-braiding $\Eps_\si$ is
given by $\lan\beta,\si\ran$. Then the connected component
$\cG_0$ of $\cG$ containing $\id\in\sys$ is the dual
principal graph of the inclusion $M\otimes M^\op\subset R$.
\labl{dualpg}
\end{proposition}

This completes our review of \cite{I1}.

\section{Quantum doubles of full induced systems}
\labl{sectMXMa}

In this section we study Longo-Rehren subfactors
$M\otimes M^\op\subset R(\sys)$ arising from the system
$\sys=\MYMa$, the full $\a$-induced system associated to
the subsystem $\NYN\subset\NXN$. In order to proceed with
$\eta$-extensions we first introduce some half-braidings.

For $\beta\in\MXMa$ choose an isometry
$T\in\Hom(\beta,\a^+_{\nu} \a^-_{\nu'})$ with
some $\nu,\nu'\in\NXN$. (These exist by definition.)
For any $\la\in\Sigma(\NXN)$ we now put
\be
\Eps^\pm_\la (\beta) = T^* \eps^\pm(\la,\nu\nu')
\a^\pm_\la (T)
\label{Eladef}
\ee
We then have the following

\begin{lemma}
The operators $\Eps^\pm_\la (\beta)$ are independent
of the choice of $T$ and $\nu,\nu'$ in the sense that,
if $\xi,\xi'\in\NXN$ and
$S\in\Hom(\beta,\a^+_{\xi}\a^-_{\xi'})$ is an isometry,
then
$\Eps^\pm_\la (\beta) = S^* \eps^\pm(\la,\xi\xi')
\a^\pm_\la (S)$. Moreover, for each $\la\in\Sigma(\NXN)$,
the family $\{\Eps^\pm_\la(\beta)\}_{\beta\in\Phi}$ is a
half-braiding for the morphism $\a^\pm_\la$ with respect
to the system $\Phi=\MXMa$.
\labl{half2}
\end{lemma}

\begin{proof}
Note that if $\beta\in\MXMa$, $\nu,\nu'\in\NXN$ and
$T\in\Hom(\beta,\a^+_{\nu}\a^-_{\nu'})$ is an
isometry, then $TT^*\in\Hom(\iota\nu\nu',\iota\nu\nu')$
since $\a^+_{\nu}\a^-_{\nu'}\iota=\iota\nu\nu'$. Hence
$TT^*\eps^\pm(\la,\nu\nu')=\eps^\pm(\la,\nu\nu')\a^\pm_\la(TT^*)$
for any $\la\in\Sigma(\NXN)$. With this it is easy to check that
$\Eps^\pm_\la(\beta)$ is unitary. The first inclusion of
\erf{ahomincl} together with \cite[Lemma 3.24]{BE1} imply
that $\eps^\pm(\la,\nu\nu')$ is an intertwiner from
$\a^\pm_\la\a^+_{\nu}\a^-_{\nu'}$ to
$\a^+_{\nu}\a^-_{\nu'}\a^\pm_\la$. With that it is easy
to check that
$\Eps^\pm_\la(\beta)\in\Hom(\a^\pm_\la\beta,\beta\a^\pm_\la)$
(cf.\ the proof of \cite[Lemma 3.20]{BE3}).
Next, for $\beta_j\in\MXMa$, $\nu_j,\nu_j'\in\NXN$ and
$T_j\in\Hom(\beta_j,\a^+_{\nu_j}\a^-_{\nu_j'})$ isometries,
$j=1,2,3$, and $X\in\Hom(\beta_3,\beta_1\beta_2)$ one has
$\a^+_{\nu_1}\a^-_{\nu_1'}(T_2)T_1XT_3^*\in
\Hom(\iota\nu_3\nu_3',\iota\nu_1\nu_1'\nu_2\nu_2')$,
and hence we can compute
\[ \begin{array}{ll}
X \Eps^\pm_\la(\beta_3)
&=
X T_3^* \eps^\pm(\la,\nu_3\nu_3') \a^\pm_\la(T_3)
= T_1^* \a^+_{\nu_1}\a^-_{\nu_1'}(T_2^*T_2)T_1X
T_3^* \eps^\pm(\la,\nu_3\nu_3') \a^\pm_\la(T_3) \\[.4em]
&= T_1^* \a^+_{\nu_1}\a^-_{\nu_1'}(T_2)^*
\nu_1\nu_1'(\eps^\pm(\la,\nu_2\nu_2'))\eps^\pm(\la,\nu_1\nu_1')
\a^\pm_\la(\a^+_{\nu_1}\a^-_{\nu_1'}(T_2)T_1X) \\[.4em]
&= T_1^* \a^+_{\nu_1}\a^-_{\nu_1'}(T_2)^*
\nu_1\nu_1'(\eps^\pm(\la,\nu_2\nu_2'))\a^+_{\nu_1}\a^-_{\nu_1'}
\a^\pm_\la(T_2)\eps^\pm(\la,\nu_1\nu_1') \a^\pm_\la(T_1X) \\[.4em]
&= \beta_1(T_2^*\eps^\pm(\la,\nu_2\nu_2')\a^\pm_\la(T_2)
T_1^*\eps^\pm(\la,\nu_1\nu_1') \a^\pm_\la(T_1X) \\[.4em]
&= \beta_1(\Eps^\pm_\la(\beta_2))\Eps^\pm_\la(\beta_1)
\a^\pm_\la(X)
\end{array} \,, \]
establishing 2.\ of Definition \ref{half}. Finally,
putting $\nu_2=\nu_2=\id$ so that consequently
$\beta_2=\id$ and $T_2=\bfe$, and choosing $X=\bfe$
gives the desired invariance properties of
$\Eps_\la(\beta)$ with $\beta=\beta_1=\beta_2$.
\end{proof}

Restricting the half-braidings to
$\sys=\MYMa\subset\MXMa=\Phi$, i.e.\ putting
$\Eps^\pm_\la=\{\Eps^\pm_\la(\beta)\}_{\beta\in\MYMa}$,
we conclude that there are extensions 
$\eta(\a^\pm_\la,\Eps^\pm_\la)$ whenever $\la\in\Sigma(\NXN)$.
Note that
\[ \Eps^\pm_\la(\beta)\a^\pm_\la(\Eps^\pm_\mu(\beta))
= T^* \eps^\pm(\la,\nu\nu')\a^\pm_\la (TT^*)
\la(\eps^\pm(\mu,\nu\nu') \a^\pm_{\la\mu} (T))
= \Eps^\pm_{\la\mu}(\beta) \]
for all $\beta\in\MXMa$, and consequently
\be
\eta(\a^\pm_{\la\mu},\Eps^\pm_{\la\mu}))=
\eta(\a^\pm_\la,\Eps^\pm_\la))
\eta(\a^\pm_\mu,\Eps^\pm_\mu))
\label{multipl}
\ee
for all $\la,\mu\in\Sigma(\NXN)$.

We now state an inclusion of intertwiner spaces which
is similar to the first inclusion in \erf{ahomincl}.

\begin{lemma}
We have
\be
\Hom(\la,\mu)\otimes \bfe  \subset 
\Hom(\eta(\a^\pm_\la,\Eps^\pm_\la),
\eta(\a^\pm_\mu,\Eps^\pm_\mu))
\ee
for any $\la,\mu\in\Sigma(\NXN)$.
\labl{inside}
\end{lemma}

\begin{proof}
Thanks to Izumi's result, Theorem \ref{irred},
and due to the first inclusion in \erf{ahomincl},
all what we have to verify is the relation
$\beta(x)\Eps^\pm_\la(\beta)=\Eps^\pm_\mu(\beta)x$
for all $\beta\in\MYMa$ whenever $x\in\Hom(\la,\mu)$.
For $\beta\in\MYMa$ there is some isometry
$T\in\Hom(\beta,\a^+_{\nu}\a^-_{\nu'})$ with some
$\nu,\nu'\in\NYN$. Then this is just
\[ \begin{array}{ll}
\beta(x)\Eps^\pm_\la(\beta) &=
\beta(x) T^*\eps^\pm(\la,\nu\nu') \a^\pm_\la (T)
= T^* \nu\nu'(x) \eps^\pm(\la,\nu\nu')
\a^\pm_\la (T) \\[.4em]
&= T^* \eps^\pm(\mu,\nu\nu') x \a^\pm_\la (T)
= T^* \eps^\pm(\mu,\nu\nu') \a^\pm_\la (T) x
=\Eps^\pm_\mu(\beta)x \,,
\end{array} \]
thanks to naturality.
\end{proof}

Immediately we obtain the following

\begin{corollary}
The map $\la\mapsto\eta(\a^\pm_\la,\Eps^\pm_\la)$,
$\la\in\Sigma(\NXN)$, preserves sums, products, and
conjugate sectors.
\labl{lahomo}
\end{corollary}

Recall from \cite[Sect.\ 4]{BE4} that
$\Hom(\id,\a^\pm_\rho)=\{ w_\rho^*v : w_\rho\in\fh_\rho^\pm \}$
where $\fh_\rho^\pm\subset\Hom(\rho,\canr)$ is the Hilbert (sub-) space
\[ \fh_\rho^\pm = \{ w_\rho\in\Hom(\rho,\canr) :
w_\rho^* \can(v) = w_\rho^* \eps^\mp (\canr,\canr) \can (v) \} \]
for $\rho\in\NXN$. Note that by Assumption \ref{assump}
the spaces $\Hom(\rho,\canr)$ and in turn $\fh_\rho^\pm$
can only be non-zero if $\rho\in\NYN$.
For any $\la,\mu\in\NXN$, $\rho\in\NYN$, we may
choose orthonormal basis of isometries
$t(_{\rho,\la}^\mu)_i \in\Hom(\mu,\rho\la)$,
$i=1,2,...,N_{\rho,\la}^\mu$ and
$w_{\rho,r;\pm}\in\fh_\rho^\pm$, where
$r=1,2,...,Z_{\rho,0}=\lan\id,\a^+_\rho\ran$
respectively $r=1,2,...,Z_{0,\rho}=\lan\id,\a^-_\rho\ran$.

\begin{lemma}
A basis of $\Hom(\a^\pm_\la,\a^\pm_\mu)$ is given by
\be
\{ t(_{\rho,\la}^\mu)_i^* w_{\rho,r;\pm}^* v : \,\,
\rho\in\NYN \,,\,\,\, i=1,2,\ldots,N_{\rho,\la}^\mu \,,\,\,\,
r=1,2,\ldots,\lan\id,\a^\pm_\rho\ran \}
\ee
for any $\la,\mu\in\NXN$.
\labl{basisaa}
\end{lemma}

\begin{proof}
It follows from $w_{\rho,r;\pm}\in\fh_\rho^\pm$ and
the first inclusion in \erf{ahomincl} that
$t(_{\rho,\la}^\mu)_i^* w_{\rho,r;\pm}^* v
\in\Hom(\a^\pm_\la,\a^\pm_\mu)$. 
The elements are clearly linearly independent as
$t(_{\rho,\la}^\mu)_i^* w_{\rho,r;\pm}^*$ are orthonormal
isometries in $N$. Now the statement follows since
$\lan\a^\pm_\la,\a^\pm_\mu\ran=\sum_\rho N_{\rho,\la}^\mu
\lan\id,\a^\pm_\rho\ran$ by Frobenius reciprocity.
\end{proof}

Next we define a subspace
$\cL(\la,\mu)\subset\Hom(\a^\pm_\la,\a^\pm_\mu)$
by putting
\be
\cL(\la,\mu)=\span \{ t(_{\rho,\la}^\mu)_i^* w_{\rho,r;\pm}^* v
: \,\,\rho\in\NYNd \,,\,\,\, i=1,2,...,N_{\rho,\la}^\mu
\,,\,\,\, r=1,2,...,Z_{\rho,0} \}
\label{cL}
\ee

Note that there is no distinction between ``$+$'' and
``$-$'' anymore because $\a^+_\rho=\a^-_\rho$
whenever $\rho\in\NYNd$.

\begin{lemma}
We have
\be
\Hom(\eta(\a^\pm_\la,\Eps^\pm_\la),\eta(\a^\pm_\mu,\Eps^\pm_\mu))
=\cL(\la,\mu) \otimes \bfe
\label{hometext}
\ee
and consequently
$\lan\eta(\a^\pm_\la,\Eps^\pm_\la),
\eta(\a^\pm_\mu,\Eps^\pm_\mu)\ran
=\sum_{\rho\in\NYNd}N_{\rho,\la}^\mu Z_{\rho,0}$
for all $\la,\mu\in\NXN$.
\labl{etalpmp}
\end{lemma}

\begin{proof}
By Theorem \ref{irred} we have to show that
\[ \cL(\la,\mu)=\{X\in \Hom(\a^\pm_\la,\a^\pm_\mu):\,\,
\Eps^\pm_\mu(\beta)X=\beta(X)\Eps^\pm_\la(\beta)
{\rm \ for\ all\ }\beta\in\MYMa \} \,. \]
So first we assume that $X$ is in the right-hand side,
and such an $X\in \Hom(\a^\pm_\la,\a^\pm_\mu)$ satisfies
$\Eps^\pm_\mu(\beta)X=\beta(X)\Eps^\pm_\la(\beta)$ in
particular for all $\beta\in\MYMpm\subset\MYMa$.
So choose an isometry
$T\in\Hom(\beta,\a^\pm_\nu)$ with some $\nu\in\NYN$.
Then, by \erf{anat},
\[ \Eps^\pm_\mu(\beta)X =T^*\eps^\pm(\mu,\nu)\a^\pm_\mu(T)X
= T^*\eps^\pm(\mu,\nu)X\a^\pm_\la(T)
= T^*\a^\mp_\nu(X) \eps^\pm(\la,\nu)\a^\pm_\la(T) \]
whereas
\[ \beta(X)\Eps^\pm_\la(\beta) =\beta(X)
T^*\eps^\pm(\la,\nu)\a^\pm_\la(T)
=T^*\a^\pm_\nu (X) \eps^\pm(\la,\nu)\a^\pm_\la(T) \,.\]
Equating these and multiplying by $T$ from the left and
$T^*$ from the right we obtain, using again \erf{anat},
\[ TT^* \a^\mp_\nu(X) \eps^\pm(\la,\nu) =
TT^* \a^\pm_\nu(X) \eps^\pm(\la,\nu) \,. \]
Since this is supposed to hold for any $\beta\in\MYMpm$
we may now take the sum over full orthonormal bases
of $\Hom(\beta,\a^\pm_\nu)$ so that we find
$\a^-_\nu(X)=\a^+_\nu(X)$ for all $\nu\in\NYN$.
Now recall that $X\in\Hom(\a^\pm_\la,\a^\pm_\mu)$
is a linear combination
\[ X=\sum_{\rho\in\NYN} \sum_{i=1}^{N_{\rho,\la}^\mu}
\sum_{r=1}^{\lan\id,\a^\pm_\rho\ran} \zeta_{\rho,i,r}
t(_{\rho,\la}^\mu)_i^* w_{\rho,r;\pm}^* v \]
with $\zeta_{\rho,i,r}\in\bbC$. But
\[ \a^\pm_\nu(t(_{\rho,\la}^\mu)_i^* w_{\rho,r;\pm}^* v)
= \nu(t(_{\rho,\la}^\mu)_i^* w_{\rho,r;\pm}^*)
\eps^\mp(\canr,\nu)v = \nu(t(_{\rho,\la}^\mu)_i^*)
\eps^\mp(\rho,\nu)w_{\rho,r;\pm}^* v \,. \]
Therefore, using $nv=0$ implies $n=0$ as well
as orthonormality of the $w_{\rho,r;\pm}$'s,
we find that $\a^-_\nu(X)=\a^+_\nu(X)$ for all $\nu$ implies
\[ \sum_{i=1}^{N_{\rho,\la}^\mu} \zeta_{\rho,i,r}
\nu(t(_{\rho,\la}^\mu)_i^*)
(\eps^+(\rho,\nu)\eps^+(\nu,\rho) - \bfe)=0 \]
for all $\nu,\rho\in\NYN$ and all
$r=1,...,\lan\id,\a^\pm_\rho\ran$.
Taking the adjoint and applying the left inverse
$\phi_\nu$ yields
\[ \sum_{i=1}^{N_{\rho,\la}^\mu}
\zeta_{\rho,i,r}^* \left( \frac{Y_{\nu,\rho}}{d_\nu d_\rho}
 -1 \right) t(_{\rho,\la}^\mu)_i=0 \,, \qquad \nu,\rho\in\NYN
\,, \quad r=1,...,\lan\id,\a^\pm_\rho\ran\,, \]
as the monodromy matrix $Y$ is obtained \cite{R1,FG,FRS2} from
$d_\nu d_\rho\phi_\nu(\eps^+(\rho,\nu)
\eps^+(\nu,\rho))^*=Y_{\nu,\rho}\bfe$.
Hence we have
$\zeta_{\rho,i,r}^* (Y_{\nu,\rho}-d_\nu d_\rho)=0$
for all $\nu,\rho,i,r$.
But $Y_{\nu,\rho}=d_\nu d_\rho$ for all $\nu\in\NYN$ 
if and only if $\rho\in\NYNd$ by \cite{R1}.
Consequently $\zeta_{\rho,i,r}=0$ whenever
$\rho\notin\NYNd$, so that indeed $X\in\cL(\la,\mu)$.

Conversely, if we start with $X\in\cL(\la,\mu)$, i.e.
\[ X=\sum_{\rho\in\NYNd} \sum_{i=1}^{N_{\rho,\la}^\mu}
\sum_{r=1}^{\lan\id,\a^\pm_\rho\ran} \zeta_{\rho,i,r}
t(_{\rho,\la}^\mu)_i^* w_{\rho,r;\pm}^* v \,,
\qquad \zeta_{\rho,i,r} \in\bbC \,, \]
then we find
$\a^-_\nu(X)=\a^+_\nu(X)$ for all $\nu\in\NYN$.
Hence, if $\beta\in\MYMa$ and
$T\in\Hom(\beta,\a^+_{\nu}\a^-_{\nu'})$ is an isometry
with some $\nu,\nu'\in\NYN$, then
\[ \begin{array}{ll}
\Eps^\pm_\mu (\beta) X &= T^* \eps^\pm(\mu,\nu\nu')
\a^\pm_\mu(T)X=T^* \eps^\pm(\mu,\nu\nu') X
\a^\pm_\la(T) \\[.4em]
&=T^* \a^\mp_{\nu\nu'}(X) \eps^\pm(\la,\nu\nu')
\a^\pm_\la(T) 
=T^* \a^+_{\nu} \a^-_{\nu'}(X) \eps^\pm(\la,\nu\nu')
\a^\pm_\la(T) \\[.4em]
&=\beta(X) T^* \eps^\pm(\la,\nu\nu') \a^\pm_\la(T)
= \beta(X) \Eps^\pm_\la(\beta)
\end{array} \]
by \erf{anat}. Thus $X$ satisfies the desired
intertwining relations.
\end{proof}

Next we compare $\eta$-extensions with different signature.

\begin{lemma}
We have
\be
\Hom(\eta(\a^+_{\la},\Eps^+_\la),\eta(\a^-_\mu,\Eps^-_\mu))
= \left\{ \begin{array}{c@{\qquad:\qquad}l}
\cL(\la,\mu)\otimes\bfe & \la,\mu\in\NYNper \\[.4em]
\{0\} & \mbox{otherwise} \end{array} \right.
\label{mixing}
\ee
for all $\la,\mu\in\NXN$.
\labl{etalpmm}
\end{lemma}

\begin{proof}
Again by Theorem \ref{irred},
we only need to show that for $\la,\mu\in\NXN$
the linear space of intertwiners
$X\in \Hom(\a^+_\la,\a^-_\mu)$ satisfying
$\Eps^-_\mu(\beta)X=\beta(X)\Eps^+_\la(\beta)$
for all $\beta\in\MYMa$ is given by $\cL(\la,\mu)$
whenever $\la,\mu\in\NYNper$ and vanishes otherwise.
Thus suppose that $X\in \Hom(\a^\pm_\la,\a^\mp_\mu)$
satisfies $\Eps^\mp_\mu(\beta)X=\beta(X)\Eps^\pm_\la(\beta)$
for all $\beta\in\MYMa$. Then in particular
\[ T^*\eps^\mp(\mu,\nu)\a^\mp_\mu(T)X=\beta(X)
T^*\eps^\pm(\la,\nu)\a^\pm_\la(T) \]
whenever $T\in\Hom(\beta,\a^\mp_\nu)$ is an isometry
and $\nu\in\NYN$. Sandwiching this with $T$ and
$\a^\pm_\la(T)^*$ yields by use of \erf{anat}
\[ \eps^\mp (\mu,\nu) \a^\mp_\mu(TT^*) X= \a^\mp_\nu(X)
\eps^\pm (\la,\mu) \a^\pm_\la(TT^*) \,, \]
and since this is supposed for any subsector
$[\beta]$ of any $[\a^\mp_\nu]$ we can sum over
orthonormal bases of $\Hom(\beta,\a^\mp_\nu)$
so that we arrive at
\[ \eps^\mp(\mu,\nu)X= \a^\mp_\nu(X)\eps^\pm(\la,\nu)
= \eps^\pm(\mu,\nu)X \qquad \mbox{for all} \,\,\,
\nu\in\NYN \,. \]
If $X\neq 0$ then $X=t^*v$ with $t\in\Hom(\mu,\canr\la)$
some necessarily non-zero multiple of an isometry.
Therefore we have found that
$\eps^+(\mu,\nu)t^*=\eps^-(\mu,\nu)t^*$ for all
$\nu\in\NYN$ implying $\mu\in\NYNper$.
Now note that if
$X\otimes\bfe\in
\Hom(\eta(\a^+_{\la},\Eps^+_\la),\eta(\a^-_\mu,\Eps^-_\mu))$
then
$X^*\otimes\bfe\in
\Hom(\eta(\a^-_{\mu},\Eps^-_\mu),\eta(\a^+_\la,\Eps^+_\la))$
so that our calculation also yields $\la\in\NYNper$.
We conclude that the intertwiner space on the left-hand side
of \erf{mixing} is zero unless $\la,\mu\in\NYNper$.
But if $\mu\in\NYNper$, then
$\a^+_\mu=\a^-_\mu$ as well as $\Eps^+_\mu=\Eps^-_\mu$,
so that clearly
$\Hom(\eta(\a^+_{\la},\Eps^+_\la),\eta(\a^-_\mu,\Eps^-_\mu))=
\Hom(\eta(\a^+_{\la},\Eps^+_\la),\eta(\a^+_\mu,\Eps^+_\mu))$.
Then the conclusion follows from Lemma \ref{etalpmp}.
\end{proof}

For $\la\in\NXN$ conjugate half-braiding
operators are given by
\[ \bar\Eps^\pm_\la (\beta) = d_\la \bar{r}_\la^*
\a^\pm_\la(\Eps^\pm_{\co\la}(\beta)^*\beta(r_\la)) \,,
\qquad \beta\in\MXMa \,, \]
where $r_\la\in\Hom(\id,\bar\la\la)$ and
$\bar{r}_\la\in\Hom(\id,\la\bar\la)$ are the R-isometries,
i.e.\ satisfying
$\la(r_\la)^*{\bar r}_\la
=\bar\la({\bar r}_\la)^*r_\la=d_\la^{-1}\bfe$.
(Recall that these isometries also serve as R-isometries
for the $\a$-induced morphisms due to the first inclusion
in \erf{ahomincl}.)

\begin{lemma}
We have $\bar\Eps^\pm_\la(\beta)=\Eps^\pm_\la(\beta)$
for all $\beta\in\MXMa$ and all $\la\in\NXN$.
\labl{conjhb}
\end{lemma}

\begin{proof}
Let $T\in\Hom(\beta,\a^+_{\nu} \a^-_{\nu'})$ be
an isometry, $\nu,\nu'\in\NXN$. Then
\[ \begin{array}{ll}
\bar\Eps^\pm_\la(\beta) &= d_\la \bar{r}_\la^*
\a^\pm_\la(\Eps^\pm_{\co\la}(\beta)^*\beta(r_\la))
=d_\la \bar{r}_\la^*
\a^\pm_\la (\a^\pm_{\bar\la}(T)^* \eps^\pm(\bar\la,\nu\nu')^*
T \beta(r_\la)) \\[.4em]
&=d_\la T^* \bar{r}_\la^* \la( \eps^\mp(\nu\nu',\bar\la)
\nu\nu' (r_\la)) \a^\pm_\la (T)
=d_\la T^* \bar{r}_\la^* \la( \bar\la(\eps^\mp(\nu\nu',\la)^*)
r_\la) \a^\pm_\la (T) \\[.4em]
&=d_\la T^*  \eps^\pm(\la,\nu\nu') \bar{r}_\la^* \la(r_\la)
\a^\pm_\la (T) =\Eps^\pm_\la(\beta) \,,
\end{array} \]
where we used the BFE
$r_\la=\bar\la(\eps^\mp(\nu\nu',\la)
\eps^\mp(\nu\nu',\bar\la)\nu\nu'(r_\la)$.
\end{proof}

Considering only $\beta\in\MYMa$, Lemma \ref{conjhb} yields
with Proposition \ref{conj} the following

\begin{corollary}
We have
$[\eta^\op(\a^\pm_\la,\Eps^\pm_\la)]=
[\eta(\a^\pm_{\bar\la},\Eps^\pm_{\bar\la})]$
for all $\la\in\NXN$.
\labl{etaopequ}
\end{corollary}

We are now ready to state the main result of this section
in the following

\begin{theorem}
We have
\be
\lan \eta(\a^+_\la,\Eps^+_\la)\eta(\a^-_\mu,\Eps^-_\mu),
\eta(\a^+_{\la'},\Eps^+_{\la'})
\eta(\a^-_{\mu'},\Eps^-_{\mu'})\ran
= \displaystyle \sum_{\nu,\xi\in\NYNper}
\sum_{\rho\in\NYNd} N_{\bar{\la'},\la}^\nu
N_{\mu',\bar\mu}^\xi N_{\nu,\xi}^\rho Z_{\rho,0}\,,
\label{compl}
\ee
for all $\la,\la',\mu,\mu'\in\NXN$.
\labl{mainir}
\end{theorem}

\begin{proof}
Using Proposition \ref{homom} and Lemma \ref{conjhb},
we can compute
\[ \begin{array}{l}
\lan \eta(\a^+_\la,\Eps^+_\la)\eta(\a^-_\mu,\Eps^-_\mu),
\eta(\a^+_{\la'},\Eps^+_{\la'})
\eta(\a^-_{\mu'},\Eps^-_{\mu'})\ran = \\[.4em]
\qquad\qquad = \lan \eta(\a^+_{\bar{\la'}},\Eps^+_{\bar{\la'}})
\eta(\a^+_\la,\Eps^+_\la),\eta(\a^-_{\mu'},\Eps^-_{\mu'})
\eta(\a^-_{\bar\mu},\Eps^-_{\bar\mu})
\ran \\[.4em]
\qquad\qquad = \sum_{\nu,\xi\in\NXN} N_{\bar{\la'},\la}^\nu
N_{\mu',\bar\mu}^\xi \lan \eta(\a^+_\nu,\Eps^+_\nu),
\eta(\a^-_\xi,\Eps^-_\xi) \ran \,,
\end{array} \]
and now the result follows by Lemma \ref{etalpmm} and
since
$\dim\cL(\nu,\xi)=
\sum_{\rho\in\NYNd}N_{\nu,\xi}^\rho Z_{\rho,0}$.
\end{proof}

Theorem \ref{mainir} has some simple consequences
in the non-degenerate case. Let us review a bit of
category language first.
A system $\cS\subset\End(Q)$ for some type III factor
$Q$ gives a strict $C^*$-tensor category
(with conjugates, subobjects, and direct sums) in the
sense of \cite{DR,LRo}, whose objects are in $\Sigma(\cS)$.
There is a natural notion of equivalence of such
categories, and two such categories $\cC$ and $\cC'$
are equivalent \cite[Prop.\ 1.1]{HY} if and only if
there is a $C^*$-tensor functor $F:\cC\rightarrow\cC'$
such that any object in $\cC'$ is isomorphic
(unitarily equivalent) to an object in the image of $F$,
and the arrow functions
$F_{\rho,\si}:\Hom(\rho,\si)\rightarrow\Hom(F(\rho),F(\si))$
are isomorphisms for any $\rho$ and $\si$ in $\cC$.
We also have a notion of direct product for two such 
strict $C^*$-tensor categories.  That is, if we have two systems
of irreducible endomorphisms of two (type III) factors $Q$ and $R$,
we have a system of irreducible endomorphisms arising as tensor
products of pairs of irreducible endomorphisms on $Q\otimes R$.
Moreover we can pass from one system of such endomorphisms on $R$
to another ``opposite'' system on the opposite algebra $R^{\op}$
naturally.

Note that the right-hand side of \erf{compl} collapses
dramatically in case that $\NYN=\NXN$ and if the original
braiding is non-degenerate, i.e.\ if $\NYNd=\{\id\}$:
Then we are simply left with Kronecker symbols
$\delta_{\la,\la'}\delta_{\mu,\mu'}$.
As shown in \cite{BEK1}, non-degeneracy of the braiding
implies $\MXMa=\MXM$, and then the global index of
$\sys=\MXMa$ is equal to the
global index $\NXN$, $[[\sys]]=[[\NXN]]$.
Theorem \ref{mainir} and Corollary \ref{etaopequ} imply that
$\{\eta(\a^+_\la,\Eps^+_\la)
\eta^\op(\a^-_\mu,\Eps^-_\mu)\}_{\la,\mu\in\NXN}$
gives a system of irreducible $R$-$R$ morphisms.
Note that each of the morphisms in this system gives
a sector arising from $\Dsys$, and with a suitable choice
of representatives in $\Dsys$ we may assume that this
system is a subsystem of $\Dsys$.
As the statistical dimension of $\eta(\a^+_\la,\Eps^+_\la)
\eta^\op(\a^-_\mu,\Eps^-_\mu)$
is $d_{\la}d_{\mu}$, we know that its global index is
equal to $[[\NXN]]^2$.
But since $[[\Dsys]]=[[\sys]]^2$,
it follows that our system is in fact the entire $\Dsys$.
With non-degeneracy, Theorem \ref{mainir} implies that
$\lan\eta(\a^+_\la,\Eps^+_\la),\eta(\a^+_\mu,\Eps^+_\mu)\ran
=\lan\la,\mu\ran$ for any $\la,\mu\in\Sigma(\NXN)$
and consequently Lemma \ref{inside} gives equalities
\[ \Hom(\la,\mu)\otimes \bfe = \Hom(\eta(\a^+_\la,\Eps^+_\la),
\eta(\a^+_\mu,\Eps^+_\mu)) \qquad \mbox{for all} \quad
\la,\mu\in\Sigma(\NXN) \,. \]
A similar statement holds for 
$\Hom(\eta^\op(\a^-_\la,\Eps^-_\la),
\eta^\op(\a^-_\mu,\Eps^-_\mu))$ and we thus have
\[ \begin{array}{l}
\Hom(\la',\mu')\otimes \Hom(\la^{\op},\mu^{\op}) \\[.4em]
\qquad\qquad =
\Hom(\eta(\a^+_{\la'},\Eps^+_{\la'})\eta^\op(\a^-_\la,\Eps^-_\la),
\eta(\a^+_{\mu'},\Eps^+_{\mu'})\eta^\op(\a^-_\mu,\Eps^-_\mu))
\end{array} \]
for all $\la',\mu',\la,\mu\in\Sigma(\NXN)$.
Let now $\cC$ be the strict $C^*$-tensor category
arising from the direct product of $\NXN$ and $(\NXN)^\op$
and $\cC'$ be the one arising from $\cD(\sys)$.
We may now introduce a functor $F:\cC\rightarrow\cC'$ which
maps any pair $(\la',\la^\op)$ to the $R$-$R$ morphism
$\eta(\a^+_{\la'},\Eps^+_{\la'})\eta^\op(\a^-_\la,\Eps^-_\la)$,
and with arrow functions $F_{(\la',\la^\op),(\mu',\mu^\op)}$
mapping
$x\otimes y\in\Hom(\la',\mu')\otimes \Hom(\la^{\op},\mu^{\op})$
to
$x\otimes y\in\Hom(\eta(\a^+_{\la'},\Eps^+_{\la'})
\eta^\op(\a^-_\la,\Eps^-_\la),\eta(\a^+_{\mu'},\Eps^+_{\mu'})
\eta^\op(\a^-_\mu,\Eps^-_\mu))$, which is obviously a
(rather trivial) $C^*$-tensor functor.
It is similarly clear that any object in $\cC'$ is
unitarily equivalent to some object in the image of $F$,
and that the arrow functions are isomorphisms.
Therefore we have the following

\begin{corollary}
If the braiding on $\NXN$ is non-degenerate,
then the strict $C^*$-tensor category given by the
system of irreducible $R$-$R$ morphisms for the Longo-Rehren
subfactor $M\otimes M^\op \subset R$ arising from the system
$\MXM$ and that given as a direct product of those arising from
the systems $\NXN$ and $(\NXN)^\op$ are equivalent.
\labl{double}
\end{corollary}

By Izumi's result, Proposition \ref{dualpg}, we find that
the irreducible $M\otimes M^\op$-$R$ morphisms arising in
our system are labelled with $\beta\in\sys=\MXMa$ and
the multiplicity of the edges between this morphism and
$\eta(\a^+_\la,\Eps^+_\la)\eta(\a^-_{\bar\mu},\Eps^-_{\bar\mu})$
is given as $\lan\beta,\a^+_\la \a^-_{\bar\mu}\ran$,
$\la,\mu\in\NXN$. Consequently the canonical endomorphism
$\Gamma\in\End(R)$ for the subfactor
$M\otimes M^\op\subset R$ decomposes as
\[ [\Gamma] = \bigoplus_{\la,\mu\in\NXN} Z_{\la,\mu} \,\,
[\eta(\a^+_\la,\Eps^+_\la)\eta^\op(\a^-_\mu,\Eps^-_\mu)] \,, \]
as $Z_{\la,\mu}=\lan\id,\a^+_\la \a^-_{\bar\mu}\ran$
and by Corollary \ref{etaopequ}.
Using the isomorphism of the tensor categories of
Corollary \ref{double} gives another

\begin{corollary}
Assume that the braiding on $\NXN$ is non-degenerate. Then
$\bigoplus_{\la,\mu\in\NXN} Z_{\la,\mu} [\la\otimes\mu^\op]$
is the sector of a canonical endomorphism for some
subfactor of $B\subset N\otimes N^\op$.
\labl{rehren}
\end{corollary}

This is a special case of a recent result of
Rehren \cite[Cor.\ 1.6]{R7}, and our method gives a new proof
of this statement by looking at the dual of the usual
Longo-Rehren subfactor arising from $\MXMa$.
(Note that a canonical endomorphism does not
determine a subfactor uniquely.  So our construction and
Rehren's might produce non-isomorphic subfactors, while
they give the same canonical endomorphism.  We expect that
these two subfactors are related by an ``$\Eps$-twist'' in
the sense of Izumi as in the remark after \cite[Prop.\ 7.3]{I1}.)

\section{Quantum doubles of chiral systems}
\labl{sectMXMpm}

In this section we study the Longo-Rehren subfactors
arising from chiral subsystems $\sys=\MYMpm\subset\MXMpm=\Phi$.
Recall from \cite[Subsect.\ 3.3]{BE3} that for
$\beta_\pm\in\MXMpm$ the operators
\[ \Eps_\rmr (\beta_+,\beta_-) =
S^*\a^-_\mu(T)^*\eps^+(\la,\mu)\a^+_\la(S)T \]
are unitaries in $\Hom(\beta_+\beta_-,\beta_-\beta_+)$
whenever $\la,\mu\in\NXN$ and $T\in\Hom(\beta_+,\a^+_\la)$
and $S\in\Hom(\beta_-,\a^-_\mu)$ are isometries,
and they do not depend on the special choices of
$\la,\mu$ and $S,T$ which realize $\beta_\pm$. Moreover,
they constitute a ``relative braiding'' between the chiral
systems $\MXMp$ and $\MXMm$. Recall that the ambichiral
system is defined as $\MXMo=\MXMp\cap\MXMm$.
For any $\tau\in\Sigma(\MXMo)$ we now put
\[ \begin{array}{l@{\qquad \mbox{for all}\quad}l}
\Eps^-_\tau(\beta) = \Eps_\rmr(\beta,\tau)^* 
& \beta\in\MXMp \,, \\[.4em]
\Eps^+_\tau(\beta) = \Eps_\rmr(\tau,\beta) 
& \beta\in\MXMm \,.
\end{array} \]
Then the following lemma plays the role of Lemma \ref{half2}.

\begin{lemma}
For each $\tau\in\Sigma(\MXMo)$, the family
$\{\Eps^\mp_\tau(\beta)\}_{\beta\in\Phi}$ is
a half-braiding with respect to the system $\Phi=\MXMpm$.
\labl{half3}
\end{lemma}

\begin{proof}
Immediate from \cite[Prop.\ 3.12]{BE3}.
\end{proof}

The restricted half-braidings
$\Eps^\mp_\tau=\{\Eps^\mp_\tau(\beta)\}_{\beta\in\sys}$
with $\Delta=\MYMpm\subset\MXMpm=\Phi$ will provide
$\eta$-extensions. Thanks to the composition rules
of the relative braiding operators we have multiplicativity
for the $\eta$-extensions,
\be
\eta(\tau\tau',\Eps^\pm_{\tau\tau'}))=
\eta(\tau,\Eps^\pm_\tau)
\eta(\tau',\Eps^\pm_{\tau'}))
\label{multipl2}
\ee
for all $\tau,\tau'\in\Sigma(\MXMo)$.
Let us now consider such $\eta$-extensions using only
$\tau\in\MXMo$ (rather than in $\Sigma(\MXMo)$).
Then it is trivial by Theorem \ref{irred} that
\[ \Hom(\eta(\tau,\Eps^\mp_\tau),\eta(\tau',\Eps^\mp_{\tau'}))
= \delta_{\tau,\tau'} \bbC \,, \]
so that all such $\eta$-extensions are irreducible and
disjoint. Note that for $R=R(\sys)$ with $\sys=\MXMpm$
we also have extensions for $\a^\pm_\la$, $\la\in\NXN$,
using the restrictions of their half-braidings
$\{\Eps^\pm_\la(\beta)\}_{\beta\in\MXMa}$ of
Section \ref{sectMXMa} to the subsystems $\MXMpm$.
By a slight abuse of notation, we also denote them by
$\eta(\a^\pm_\la,\Eps^\pm_\la)$.
In order to avoid too confusing $\pm$-indices,
we now better focus on the case $\sys=\MXMp$.
The other case, $\sys=\MXMm$, is of course completely analogous.
The following lemma is the analogue of Lemma \ref{etalpmp},
now addressing the extensions $\eta(\a^+_\la,\Eps^+_\la)$
for $\sys=\MXMp$.

\begin{lemma}
We have
\be
\Hom(\eta(\a^+_\la,\Eps^+_\la),\eta(\a^+_\mu,\Eps^+_\mu))
=\cL(\la,\mu) \otimes \bfe
\label{hometextc}
\ee
and consequently
$\lan\eta(\a^+_\la,\Eps^+_\la),
\eta(\a^+_\mu,\Eps^+_\mu)\ran
=\sum_{\rho\in\NYNd}N_{\rho,\la}^\mu Z_{\rho,0}$
for all $\la,\mu\in\NXN$.
\labl{etalpmpc}
\end{lemma}

\begin{proof}
Literally the same as the proof of Lemma \ref{etalpmp},
apart from the simplification in the second half that we
now only have to consider $\nu'=\id$.
\end{proof}

Next we compare our different kinds of $\eta$-extensions.

\begin{lemma}
We have
\be
\Hom(\eta(\tau,\Eps^-_\tau)),\eta(\a^+_{\la},\Eps^+_\la))
= \left\{ \begin{array}{c@{\qquad:\qquad}l}
\Hom(\tau,\a^+_\la)\otimes\bfe & \la\in\NYNper \\[.4em]
\{0\} & \mbox{otherwise} \end{array} \right.
\label{mixing2}
\ee
for all $\tau\in\MXMo$ and all $\la\in\NXN$.
\labl{etalptm}
\end{lemma}

\begin{proof}
Using once again Theorem \ref{irred}, we first
show that if a non-zero $X\in\Hom(\tau,\a^+_\la)$
satisfies $\Eps^+_\la(\beta)X=\beta(X)\Eps^-_\tau(\beta)$
for all $\beta\in\MYMp$, then this implies $\la\in\NYNper$.
So suppose we have such an $X\neq0$. Since $\tau\in\MXMo$
there will also be some $\mu\in\NXN$ and an isometry
$Q\in\Hom(\tau,\a^-_\mu)$. Then the intertwining
condition reads
\[ T^*\eps^+(\la,\nu)\a^+_\la(T)X=\beta(X)
T^*\a^+_\nu(Q)^*\eps^-(\mu,\nu)\a^-_\mu(T)Q \]
whenever $\nu\in\NYN$ and $T\in\Hom(\beta,\a^+_\nu)$
is an isometry. Multiplication with $T$ from the
right yields
\[ \eps^+(\la,\nu)\a^+_\la(T)X
=\a^+_\nu(XQ^*)\eps^-(\mu,\nu)\a^-_\mu(T)Q
=\eps^-(\la,\nu)\a^+_\la(T)X \,, \]
where we exploited
$XQ^*\in\Hom(\a^-_\mu,\a^+_\la)\subset\Hom(\iota\mu,\iota\la)$
to apply \erf{anat}. Now we can multiply by $\tau(T)^*$ from
the right, and then we may use a summation over full
orthonormal bases of $\Hom(\beta,\a^+_\nu)$ to obtain
$\eps^+(\la,\nu)X=\eps^-(\la,\nu)X$ for all $\nu\in\NYN$.
Now note that a non-zero $X\in\Hom(\tau,\a^+_\la)$ is
necessarily of the from $X=t^*v$ with
$t\in\Hom(\la,\bar\iota\tau\iota)$ a non-zero multiple
of an isometry. Hence we find
$\eps^+(\la,\nu)=\eps^-(\la,\nu)$ for all $\nu\in\NYN$,
proving that the left-hand side of \erf{mixing2} is zero
unless $\la\in\NYNper$.

On the other hand, if $\la\in\NYNper$ then $\a^+_\la=\a^-_\la$.
Hence, for an arbitrary $X\in\Hom(\tau,\a^+_\la)$ we
find $X^*\in\Hom(\a^-_\la,\tau)$, and therefore the naturality
of the relative braiding of \cite[Prop.\ 3.12]{BE3}
gives us
$\Eps_\rmr(\beta,\tau)\beta(X)^*=X^*\Eps_\rmr(\beta,\a^-_\la)$
for any $\beta\in\MXMp$. By taking adjoints this reads
\[ \beta(X)\Eps^-_\tau(\beta)=T^*\eps^-(\la,\nu)\a^-_\la(T)X
=T^*\eps^+(\la,\nu)\a^+_\la(T)X=\Eps^+_\la(\beta)X \,,\]
so that the desired intertwining relation is automatically
fulfilled in particular for $\beta\in\MYMp$.
This completes the proof.
\end{proof}

Conjugate half-braiding operators are given for $\tau\in\MXMo$ by
\[ \bar\Eps^-_\tau (\beta) = d_\tau \bar{R}_\tau^*
\tau(\Eps^-_{\co\tau}(\beta)^*\beta(R_\tau)) \,,
\qquad \beta\in\MXMp \,, \]
with R-isometries $R_\tau\in\Hom(\id,\bar\tau\tau)$ and
$\bar{R}_\tau\in\Hom(\id,\tau\bar\tau)$ satisfying
$\tau(R_\tau)^*\bar{R}_\tau
=\bar\tau(\bar{R}_\tau)^*R_\tau=d_\tau^{-1}\bfe$.

\begin{lemma}
We have $\bar\Eps^-_\tau(\beta)=\Eps^-_\tau(\beta)$
for all $\beta\in\MXMp$ and all $\tau\in\MXMo$.
\labl{conjhbc}
\end{lemma}

\begin{proof}
We compute
\[ \bar\Eps^-_\tau(\beta) = d_\tau \bar{R}_\tau^*
\tau(\Eps_\rmr(\beta,{\bar\tau})\beta(R_\tau))
= d_\tau \bar{R}_\tau^*
\tau(\bar\tau(\Eps_\rmr(\beta,\tau)^* R_\tau)
= \Eps^-_\tau(\beta) \,, \]
where we used the BFE for the relative braiding
\cite[Prop.\ 3.12]{BE3},
$R_\tau=\bar\tau(\Eps_\rmr(\beta,\tau)
\Eps_\rmr(\beta,\bar\tau)\beta(R_\tau)$.
\end{proof}

Considering only $\beta\in\MYMp$, Lemma \ref{conjhbc} yields
with Proposition \ref{conj} the following

\begin{corollary}
We have
$[\eta^\op(\tau,\Eps^-_\tau)]=
[\eta(\bar\tau,\Eps^-_{\bar\tau})]$
for all $\tau\in\MXMo$.
\labl{etaopequc}
\end{corollary}

Recall from \cite{BEK2} that
$b^\pm_{\tau,\la}=\lan\tau,\a^\pm_\la\ran$
denote the chiral branching coefficients for
ambichiral $\tau$ and $\la\in\NXN$.

\begin{theorem}
We have
\be
\lan \eta(\a^+_\la,\Eps^+_\la)\eta(\tau,\Eps^-_\tau),
\eta(\a^+_\mu,\Eps^+_\mu)\eta(\tau',\Eps^-_{\tau'})\ran=
\sum_{\tau''\in\MXMo} \sum_{\rho\in\NYNper} 
N_{\bar\tau,\tau'}^{\tau''}  N_{\la,\bar\mu}^\rho
b^+_{\tau'',\rho}
\label{complc}
\ee
for all $\la,\mu\in\NXN$ and all $\tau,\tau'\in\MXMo$.
\labl{mainirc}
\end{theorem}

\begin{proof}
Analogous to the proof of Theorem \ref{mainir},
this is reduced to Lemma \ref{etalptm} by use of
Proposition \ref{homom}.
\end{proof}

Let $\Ups\subset\cD(\MYMp)$ denote the subset of morphisms
$\Om\in\cD(\MYMp)$ which correspond to subsectors of
$[\eta(\a^+_\la,\Eps^+_\la)\eta(\tau,\Eps^-_\tau)]$
considering $\la\in\NYN$ and $\tau\in\MYMo$ only.
Now the question arises whether $\Ups$ is a proper
subsystem or whether it may exhaust the entire quantum
double system and therefore we would like to measure its size.
For this purpose we compare the global indices
$[[\Ups]]$ and $[[\cD(\MYMp)]]=[[\MYMp]]^2$.

\begin{proposition}
The global index of $\Ups$ is given by
\be
[[\Ups]] = \frac {\sum_{\rho\in\NYNd} d_\rho Z_{\rho,0}}
{[[\NYNd]]} \, [[\cD(\MYMp)]] \, .
\label{uuuups}
\ee
\labl{upsglob}
\end{proposition}

\proof
Let $R_{\tau,\la}$, $\tau\in\MYMo$, $\la\in\NYN$, denote matrices
with entries
\[ R_{\tau,\la;\Om}^{\Om'}=
\lan \Om\eta(\a^+_\la,\Eps^+_\la) \eta(\tau,\Eps^-_\tau),
\Om'\ran \,,\qquad \Om,\Om'\in\Ups \,.\]
Further let $\vec{d}$ denote the column vector
with entries $d_\Om$, $\Om\in\Ups$.
Then $\vec{d}$ is a simultaneous eigenvector of the
matrices $R_{\tau,\la}$ with respective eigenvalues
$d_\tau d_\la$. We define another vector $\vec{v}$
by putting
\[ v_\Om = \sum_{\tau\in\MYMo} \sum_{\la\in\NYN} d_\tau d_\la
\lan \Om, \eta(\a^+_\la,\Eps^+_\la)
\eta(\tau,\Eps^-_\tau) \ran \,, \qquad \Om\in\Ups \,. \]
Then we have $R_{\tau,\la}\vec{v}=d_\tau d_\la \vec{v}$,
as we can compute
\[ \begin{array}{ll}
(R_{\tau,\la} \vec{v})_\Om &= \sum_{\Om'\in\Ups}
\sum_{\tau'\in\MYMo} \sum_{\mu\in\NYN}
\lan \Om \eta(\a^+_\la,\Eps^+_\la)
\eta(\tau,\Eps^-_\tau) , \Om'\ran \\[.4em]
& \qquad\qquad\qquad\qquad\qquad \times d_{\tau'} d_\mu
\lan \Om', \eta(\a^+_\mu,\Eps^+_\mu)
\eta(\tau',\Eps^-_{\tau'}) \ran  \\[.4em]
&= \sum_{\tau'\in\MYMo} \sum_{\mu\in\NYN}
d_{\tau'} d_\mu \lan \Om ,
\eta(\a^+_{\bar\la} \Eps^+_{\bar\la})
\eta(\a^+_\mu,\Eps^+_\mu) \eta(\tau',\Eps^-_{\tau'}) 
\eta(\bar\tau,\Eps^-_{\bar\tau})\ran \\[.4em]
&= \sum_{\tau',\tau''\in\MYMo}
\sum_{\mu,\nu\in\NXN} d_{\tau'} d_\mu
N_{\bar\la,\mu}^\nu N_{\tau',\bar\tau}^{\tau''}
\lan \Om, \eta(\a^+_\nu,\Eps^+_\nu)
\eta(\tau'',\Eps^-_{\tau''}) \ran \\[.4em]
&= \sum_{\tau''\in\MYMo} \sum_{\nu\in\NXN}
d_\tau d_\la d_{\tau''} d_\nu
\lan \Om, \eta(\a^+_\nu,\Eps^+_\nu)
\eta(\tau'',\Eps^-_{\tau''}) \ran =
d_\tau d_\la v_\Om \,.
\end{array} \]
Because the sum matrix $\sum_{\tau,\la} R_{\tau,\la}$ is
irreducible it follows $\vec{v}=\zeta \vec{d}$,
$\zeta\in\bbR$, by the uniqueness of the Perron-Frobenius
eigenvector. Note that
\[ d_\tau d_\la =\sum_{\Om\in\Ups} 
\lan\Om, \eta(\a^+_\la,\Eps^+_\la)
\eta(\tau,\Eps^-_\tau) \ran d_\Om \,, \]
and hence
$[[\NYN]][[\MYMo]]=\sum_\Om v_\Om d_\Om = \zeta [[\Ups]]$.
We next notice that $\zeta=v_\id$ as $d_\id=1$. But $v_\id$
can be computed as
\[ v_\id = \sum_{\tau\in\MYMo} \sum_{\la\in\NYN} d_\tau d_\la 
\lan \eta(\bar\tau,\Eps^-_{\bar\tau}) ,
\eta(\a^+_\la,\Eps^+_\la) \ran =
\sum_{\la\in\NYNd} d_\la \sum_{\tau\in\MYMo}
b^+_{\tau,\la} d_\tau =\sum_{\la\in\NYNd} d_\la^2 \,, \]
where we used Lemma \ref{etalptm}.
Hence $[[\Ups]]=[[\NYN]][[\MYMo]]/[[\NYNd]]$, and now
the claim follows since
$[[\MYMo]]=(\sum_{\rho\in\NYNd} d_\rho Z_{\rho,0})
[[\MYMp]]^2/[[\NYN]]$
by \cite[Prop.\ 3.1]{BE4}.
\endproof

Similar to Theorem \ref{mainir}, the degenerate
morphisms $\rho$ appearing in \erf{complc} are
responsible that some of the
$\eta(\a^+_\la,\Eps^+_\la)\eta(\tau,\Eps^-_\tau)$'s
will be equivalent or are reducible,
and this will cause some kind of orbifolding
as we will show in Section \ref{col0} by examples.
Note, however, that also the right-hand side
of \erf{complc} simplifies considerably if the original
braiding is non-degenerate and if we have $\NXN=\NYN$:
We are just left with Kronecker symbols
$\delta_{\la,\mu}\delta_{\tau,\tau'}$.
Since the statistical dimension of $\eta(\tau,\Eps^-_\tau)$ is
$d_\tau$ and as $[[\NXN]][[\MXMo]]=[[\MXMp]]^2$ thanks
to \cite[Thm.\ 4.2]{BEK2}, we conclude that the family of morphisms
$\eta(\a^+_\la,\Eps^+_\la)\eta^\op(\tau,\Eps^-_\tau)$ serves
as a system $\Dsys$. In the non-degenerate case, it is
derived similarly from Theorem \ref{mainirc} and
Corollary \ref{etalptm} that then
\[ \Hom(\la,\mu)\otimes \bfe = \Hom(\eta(\a^+_\la,\Eps^+_\la),
\eta(\a^+_\mu,\Eps^+_\mu)) \qquad \mbox{for all} \quad
\la,\mu\in\Sigma(\NXN) \,, \]
as well as
\[ \bfe \otimes \Hom(\tau^\op,{\tau'}^\op)
= \Hom(\eta^\op(\tau,\Eps^-_\tau),
\eta^\op(\tau',\Eps^-_{\tau'})) \qquad \mbox{for all} \quad
\tau,\tau'\in\Sigma(\MXMo) \,. \]
By the same arguments which lead to Corollary \ref{double}
this gives  the following

\begin{corollary}
If the braiding on $\NXN$ is non-degenerate, then the
strict $C^*$-tensor category given by the system of
irreducible $R$-$R$ morphisms for the Longo-Rehren subfactor
$M\otimes M^\op \subset R$ arising from the system $\MXMpm$
and that given as a direct product of those arising from
the systems $\NXN$ and $(\MXMo)^\op$ are equivalent.
\labl{double2}
\end{corollary}

This corollary seems to be a precise statement of an
announcement by Ocneanu. Namely, at the Taniguchi Conference
in Nara, Japan, in December 1998, he announced as a part
of his ``big sandwich of theorems'' that ``the quantum double
of a quantum subgroup$^+$ of a non-degenerately braided
quantum group is equal to the quantum group
$\times\; \overline{\rm ambichirals}$'' (in whatever sense).

Note that the braiding on the ``quantum double''
system of $R$-$R$ morphisms is given by the direct product
of the original one on $\NXN$ and the one
on $\MXMo$ in the above theorem.
Since \cite[Thm.\ 5.5]{I1} implies that this braiding is
non-degenerate and since we assumed non-degeneracy of
the original braiding on $\NXN$ here, Corollary \ref{double2}
implies that also the braiding on the ambichiral system
$\MXMo$ is non-degenerate, in perfect agreement with
our result \cite[Thm.\ 4.2]{BEK2}.

The A-D-E cases studied in \cite{O3,X1,BE2,BEK2} provide
the following examples.

\begin{corollary}
As strict $C^*$-tensor categories,
the quantum double systems of the chiral systems
$\mathrm{E}_6$, $\mathrm{E}_8$, and $\mathrm{D}_{2n}$ are
equivalent to $\mathrm{A}_{11}\times (\mathrm{A}_3)^\op$,
$\mathrm{A}_{29}\times (\mathrm{A}_4^{\mathrm{even}})^\op$,
and
$\mathrm{A}_{4n-3}\times (\mathrm{D}_{2n}^{\mathrm{even}})^\op$,
respectively.
\labl{ADE}
\end{corollary}

By the same arguments used in Section \ref{sectMXMa}
we now find for the non-degenerate case and
$\NXN=\NYN$ that
\[ [\Gamma]=\bigoplus_{\la\in\NXN}\bigoplus_{\tau\in\MXMo}
b^+_{\tau,\la} \,\,
[\eta(\a^+_\la,\Eps^+_\la)\eta^\op(\tau,\Eps^-_\tau)] \]
is the canonical endomorphisms sector of
$M\otimes M^\op\subset R$ arising from $\sys=\MXMp$.
However, if one considers the Longo-Rehren subfactor
arising from $\sys=\MYMp$ where $\NYN$ is now a proper
and degenerate subsystem of $\NXN$, then the computations
for the structure of $\Dsys$ and the dual principal
graph become more involved. In that case one needs
the whole general machinery of this section which takes
care of possible degeneracies. Such situations will
be handled in Section \ref{col0}.

\section{Quantum doubles of color zero subsystems of
chiral systems}
\labl{col0}

Subfactors with principal graphs E$_6$, E$_8$ are
basic and important examples of subfactors
arising from $\a$-induction \cite{X1,BE2}.
The Longo-Rehren subfactor arising from the subfactor
with principal graph E$_6$ has been studied and the principal
and the dual principal graphs have been computed, as well as
other information, by Izumi \cite{I2}. Note that this
subfactor is {\sl different} from the Longo-Rehren
subfactor arising from the chiral system for the
conformal inclusion $SU(2)_{10}\subset SO(5)_1$
as studied in Section \ref{sectMXMpm}.
The reason is that Izumi considers in \cite{I2} the
quantum double system of the endomorphisms corresponding
to the {\sl even} three vertices rather than all
nodes of the graph E$_6$.
This is more natural from the viewpoint
of the usual theory of type II$_1$ subfactors,
since we obtain this quantum double system of the system
of the three $M$-$M$ bimodules, if we apply
the construction of the asymptotic inclusion
$M\vee (M'\cap M_\infty)\subset M_\infty$
to the hyperfinite II$_1$ subfactor $N\subset M$
with principal graph E$_6$ and compute the system of
$M_\infty$-$M_\infty$ bimodules.
So we will study the Longo-Rehren subfactors arising from
$\a$-induction corresponding to this type of
asymptotic inclusions in this section.
That is, from the view point of the $\a$-induction,
the chiral system $\MXMp$ for the subfactor $N\subset M$
arising from the conformal inclusion
$SU(2)_{10}\subset SO(5)_1$ has a natural ``coloring'' for
irreducible objects with colors 0 and 1, inherited from
the coloring of the $\SUz_{10}$ system coming from the
even-odd parity of the spins.
More precisely and generally, thanks to Wassermann's work
\cite{W}, we know that there are (non-degenerately) braided
systems $\cX_{n,k}=\{\la_\La:\La\in\cA_{n,k}\}$,
where $\cA_{n,k}$ denotes the $\SUn$ level $k$ Weyl alcove,
such that the morphisms $\la_\La\in\End(N)$ satisfy the
$\SUn_k$ fusion rules and have statistics phases
$\om_\La=\E^{2\pi\I h_\La}$, where $h_\La$ are the
conformal dimensions, for any $n,k=1,2,...\,$.
The Weyl alcove has a natural coloring (``$n$-ality'')
$t:\cA_{n,k}\rightarrow\bbZ_n$, and the
color zero subsystems $\cY_{n,k}\subset\cX_{n,k}$
are given by $\cY_{n,k}=\{\la_\La:t(\La)=0\}$.
Now let $N\subset M$ be a subfactor arising from a
conformal inclusion $SU(n)_k\subset G_1$ for some
Lie group $G$, as treated in \cite{BEK2}.
We put $\NXN=\cX_{n,k}$.
Note that then the ambichiral system $\MXMo$ corresponds
to the positive energy representations $\pi_\ell$ of $G_1$,
and the chiral branching coefficients are the
well-known branching coefficients of the conformal
inclusion at hand
$b^\pm_{\tau_\ell,\la_\La}=b_{\ell,\La}$
and the modular invariant matrix is given by
\[ Z_{\La,\La'}=\sum_\ell b_{\ell,\La} b_{\ell,\La'} \,.\]
We now set $\Phi=\MXMp$ and $\sys=\MYMp$, where $\MYMp$
arises from the color zero subsystem $\NYN=\cY_{n,k}$.
Note that $\NYN$ will in general be degenerate though
we have always non-degeneracy for $\NXN$ here.
Here we will study Longo-Rehren subfactors
$M\otimes M^\op\subset R$ with $R=R(\sys)$
and illustrate that the degeneracy of $\NYN$ causes
naturally a certain orbifold procedure by means of
Theorem \ref{mainirc}.
The examples we cover correspond to subfactors with
principal graph E$_6$ and E$_8$, and all the three analogues
arising from conformal inclusions of $SU(3)_k$.

\begin{example}\labl{E6}{\rm
We start with the subfactor $N\subset M$
arising from the conformal inclusion
$SU(2)_{10}\subset SO(5)_1$.
The irreducible endomorphisms in $\la_j\in\NXN$ are
labelled with $j\in\{0,1,\dots,10\}$ as usual and those
in $\MXMo$ are labelled with $\tau_\ell$, $\ell=0,1,2$
as the vertices of A$_3$. (Such that $\tau_0=\id$.)
The morphisms $\tau_\ell$ obey Ising fusion rules,
$[\tau_1\tau_1]=[\tau_0]\oplus[\tau_2]$, the non-vanishing
branching coefficients $b^+_{\tau_\ell,\la_j}=b_{\ell,j}$
are given by
\[ b_{0,0}=b_{0,6}=b_{1,3}=b_{1,7}=b_{2,4}=b_{2,10}=1 \,,\]
and the E$_6$ modular matrix is given by
$Z_{j,j'}=\sum_{\ell=0}^2 b_{\ell,j} b_{\ell,j'}$
(cf.\ \cite[Example 2.2]{BE2}).
With $\NYN=\{\la_j:j=0,2,4,6,8,10\}$,
we study the Longo-Rehren subfactor
$M\otimes M^\op\subset R$ arising from the
system $\sys=\MYMp$ whose irreducible morphisms correspond
to the three even vertices of the graph E$_6$.
Note that $\la_{10}$ is degenerate in $\NYN$, in fact
we have $\NYNd=\NYNper=\{\la_0,\la_{10}\}$.
Since $Z_{10,0}=0$ we find by Proposition \ref{upsglob}
that the set $\Ups$ is only provides half of the
quantum double system $\cD(\MYMp)$.
Thus, as considering even $j$ and even $\ell$ only will not
exhaust $\cD(\MYMp)$, we now consider $\eta(\a^+_j,\Eps^+_j)$,
$\eta(\tau_\ell,\Eps^-_\ell)$ with
$j\in\{0,1,\dots,10\}$ and $\ell\in\{0,1,2\}$,
where we write $\a_j$ for $\a_{\la_j}$ as
in \cite[Example 2.2]{BE2}.
These extended endomorphisms of $R$ may not decompose
into direct sums of irreducible morphisms in $\Dsys$ any more,
but $\eta(\a^+_j,\Eps^+_j)\eta(\tau_\ell,\Eps^-_\ell)$ do
decompose into direct sums of irreducible morphisms in $\Dsys$
if $j+\ell$ is even, since then $\a^+_j\tau_\ell$ decompose
into a direct sum of morphisms in $\MYMp$.
Now Lemma \ref{etalpmpc} yields easily irreducibility
of $\eta(\a^+_j,\Eps^+_j)$ for $j\in\{0,2,4,6,8,10\}$.
Lemma \ref{etalptm} yields similarly
\[ \lan \eta(\a^+_{10},\Eps^+_{10}),
\eta(\tau_2,\Eps^-_2)\ran=1 \,,\]
and by Theorem \ref{mainirc} we find more generally
\[ \lan \eta(\a^+_j,\Eps^+_j)\eta(\tau_0,\Eps^-_0),
\eta(\a^+_{10-j},\Eps^+_{10-j})\eta(\tau_2,\Eps^-_2)\ran = 1 \]
for $j=0,2,4,...,10$. We similarly compute
\begin{eqnarray*}
\lan \eta(\a^+_1,\Eps^+_1)\eta(\tau_1,\Eps^-_1),
\eta(\a^+_1,\Eps^+_1)\eta(\tau_1,\Eps^-_1)\ran &=& 1,\\
\lan \eta(\a^+_3,\Eps^+_3)\eta(\tau_1,\Eps^-_1),
\eta(\a^+_3,\Eps^+_3)\eta(\tau_1,\Eps^-_1)\ran &=& 1,\\
\lan \eta(\a^+_5,\Eps^+_5)\eta(\tau_1,\Eps^-_1),
\eta(\a^+_5,\Eps^+_5)\eta(\tau_1,\Eps^-_1)\ran &=& 2,\\
\lan \eta(\a^+_7,\Eps^+_7)\eta(\tau_1,\Eps^-_1),
\eta(\a^+_7,\Eps^+_7)\eta(\tau_1,\Eps^-_1)\ran &=& 1,\\
\lan \eta(\a^+_9,\Eps^+_9)\eta(\tau_1,\Eps^-_1),
\eta(\a^+_9,\Eps^+_9)\eta(\tau_1,\Eps^-_1)\ran &=& 1,\\
\lan \eta(\a^+_1,\Eps^+_1)\eta(\tau_1,\Eps^-_1),
\eta(\a^+_9,\Eps^+_9)\eta(\tau_1,\Eps^-_1)\ran &=& 1,\\
\lan \eta(\a^+_3,\Eps^+_3)\eta(\tau_1,\Eps^-_1),
\eta(\a^+_7,\Eps^+_7)\eta(\tau_1,\Eps^-_1)\ran &=& 1.
\end{eqnarray*}
We have a decomposition of
$\eta(\a^+_5,\Eps^+_5)\eta(\tau_1,\Eps^-_1)$
into two irreducible, mutually inequivalent endomorphisms
$\Omega,\Omega'\in\End(R)$ which must
belong (up to equivalence) to $\Dsys$.
By Izumi's result \cite[Lemma 4.5]{I1} we conclude
that that there are morphisms $\si,\si'\in\Sigma(\MYMp)$
with half-braidings $\Eps_\si,\Eps_{\si'}'$ such
that $\Omega=\eta(\si,\Eps_\si)$ and
$\Omega'=\eta(\si',\Eps_{\si'}')$.
But since we have
$[\alpha_5^+][\tau_1]=2[\alpha_2^+]$
and since $\eta$-extension preserves the statistical
dimension we must have
$d_\Omega+d_{\Omega'}=d_5 d_{\tau_1}=2d_2$.
(Recall $d_{\a^+_j}=d_j$.)
Moreover, Theorem \ref{irred} implies that
$[\a^+_2]$ is a subsector of both, $[\si]$ and $[\si']$,
and in turn $d_\Omega=d_\sigma\ge d_2$,
$d_{\Omega'}=d_{\sigma'}\ge d_2$.
This forces $d_\Omega=d_{\Omega'}=d_2$ and consequently
$[\a^+_2]=[\si]=[\si']$, i.e.\ the product
$\eta(\a^+_5,\Eps^+_5)\eta(\tau_1,\Eps^-_1)$ decomposes
into two irreducibles of equal statistical dimension.
So we have (at least) the following 8 irreducible,
mutually inequivalent endomorphisms
$\eta(\a^+_0,\Eps^+_0)$,
$\eta(\a^+_2,\Eps^+_2)$,
$\eta(\a^+_4,\Eps^+_4)$,
$\eta(\a^+_6,\Eps^+_6)$,
$\eta(\a^+_8,\Eps^+_8)$,
$\eta(\a^+_{10},\Eps^+_{10})$,
$\eta(\a^+_1,\Eps^+_1)\eta(\tau_1,\Eps^-_1)$,
$\eta(\a^+_3,\Eps^+_3)\eta(\tau_1,\Eps^-_1)$
and two more irreducible endomorphisms
$\Omega,\Omega'$ of $R$ arising from
$\eta(\a^+_5,\Eps^+_5)\eta(\tau_1,\Eps^-_1)$. 
Counting the global index, we conclude that these are all the 
$R$-$R$ morphisms in the system $\Dsys$.
Thus, with Proposition \ref{dualpg} and recalling that
$\eta(\a^+_j,\Eps^+_j)\eta(\tau_\ell,\Eps^-_\ell)$
is the $\eta$-extension of $\a^+_j\circ\tau_\ell$
with the composed half-braiding, we can compute the dual
principal graph of the subfactor $M\otimes M^\op\subset R$
which we display in Fig.\ \ref{dE6}.
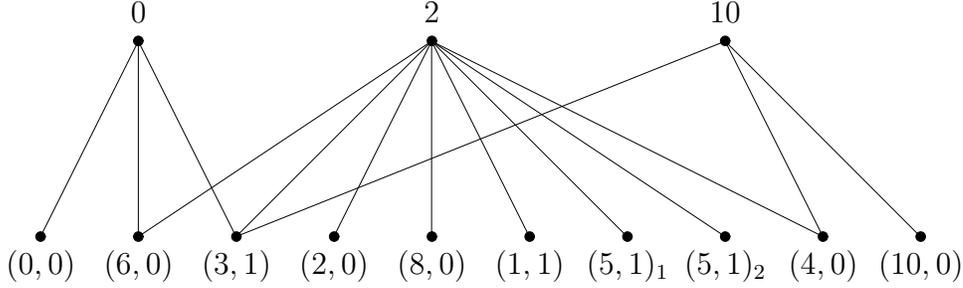
\begin{figure}[htb]
\begin{center}
\unitlength 1.3mm
\begin{picture}(110,40)
\thinlines 
\multiput(10,10)(10,0){10}{\circle*{1}}
\multiput(20,30)(30,0){3}{\circle*{1}}
\path(20,30)(10,10)
\path(20,30)(20,10)
\path(20,30)(30,10)
\path(50,30)(20,10)
\path(50,30)(30,10)
\path(50,30)(40,10)
\path(50,30)(50,10)
\path(50,30)(60,10)
\path(50,30)(70,10)
\path(50,30)(80,10)
\path(50,30)(90,10)
\path(80,30)(30,10)
\path(80,30)(90,10)
\path(80,30)(100,10)
\put(10,7){\makebox(0,0){$(0,0)$}}
\put(20,7){\makebox(0,0){$(6,0)$}}
\put(30,7){\makebox(0,0){$(3,1)$}}
\put(40,7){\makebox(0,0){$(2,0)$}}
\put(50,7){\makebox(0,0){$(8,0)$}}
\put(60,7){\makebox(0,0){$(1,1)$}}
\put(70,7){\makebox(0,0){$(5,1)_1$}}
\put(80,7){\makebox(0,0){$(5,1)_2$}}
\put(90,7){\makebox(0,0){$(4,0)$}}
\put(100,7){\makebox(0,0){$(10,0)$}}
\put(20,33){\makebox(0,0){$0$}}
\put(50,33){\makebox(0,0){$2$}}
\put(80,33){\makebox(0,0){$10$}}
\end{picture}
\end{center}
\caption{The dual principal graph for the
Longo-Rehren subfactor arising from E$_6$}
\label{dE6}
\end{figure}
Of course it is the same as the one first
computed by Izumi \cite{I2} by direct computations
of the tube algebra involving $6j$-symbols.
In the graph, we used an obvious notation for the
vertices labelled by morphisms in $\MYMp$, and we
simply wrote $(j,\ell)$ for the pair
$\eta(\a^+_j,\Eps^+_j)\eta(\tau_\ell,\Eps^-_\ell)$ with
$j\in\{0,1,\dots,10\}$ and $\ell\in\{0,1,2\}$. The labels
$(5,1)_1$ and $(5,1)_2$ stand for the two irreducible
endomorphisms $\Omega,\Omega'$ corresponding to the
subsectors of $\eta(\a^+_5,\Eps^+_5)\eta(\tau_1,\Eps^-_1)$.
The procedure yielding the $R$-$R$ morphisms
here is an orbifold procedure of order 2 for the
$(j,\ell)$ with $j+\ell\in 2\bbN$ with symmetry
$(j,\ell)\leftrightarrow (10-j,2-\ell)$.}
\end{example}

\begin{example}\labl{E8}{\rm 
We next study the Longo-Rehren subfactor arising from the
four even vertices of the graph E$_8$ and compute the 
dual principal graph, which is new.
The subfactor $N\subset M$ now arises from the
conformal inclusion $SU(2)_{28}\subset (\Gtwo)_1$
as in \cite[Example 2.3]{BE2}. Analogously to
Example \ref{E6}, the full system is given
$\NXN=\{\la_j\mid j=0,1,2,...,28\}$.
We label the ambichiral morphisms in $\MXMo$ with
$\tau_k$, $k=0,2$, corresponding to the extremal
vertices of the two long legs of E$_8$.
They $\tau_\ell$'s obey Lee-Yang fusion rules,
$[\tau_2\tau_2]=[\tau_0]\oplus[\tau_2]$ which
is the fusion of the even vertices of A$_4$.
The non-vanishing branching coefficients
$b^+_{\tau_\ell,\la_j}=b_{\ell,j}$
are given by
\[ b_{0,0}=b_{0,10}=b_{0,18}=b_{0,28}=
b_{2,6}=b_{2,12}=b_{2,16}=b_{2,22}=1 \,,\]
determining the modular invariant $Z$ as before.
The color zero subsystem is
$\NYN=\{\la_j\mid j=0,2,4,...,28\}$, and then
$\NYNd=\{\la_0,\la_{28}\}$.
We will study the Longo-Rehren subfactor
$M\otimes M^\op\subset R$ arising from the chiral
induced system $\sys=\MYMp$.
The system $\Phi=\MXMp$ corresponds to the labels of
vertices of E$_8$ in \cite[Fig.\ 8]{BE2}, and the
subsystem $\sys=\MYMp$ to the even ones.
As the ambichiral vertices are both even we find
$\MXMo=\MYMo$ here.
Note that the degenerate morphism $\la_{28}$
appears in the vacuum column and row of $Z$ this time.
Therefore Proposition \ref{upsglob} tells us that
$\Ups=\Dsys$, i.e.\ in contrast to Example \ref{E6}
we do not need to consider the odd spins at all in order
to produce the entire $\Dsys$ by our method.
Lemma \ref{etalpmpc} gives
$[\eta(\a^+_j,\Eps^+_j)]=[\eta(\a^+_{28-j},\Eps^+_{28-j})]$
for $j=0,2,4,\dots,12$ and it similarly implies that
$\eta(\a^+_j,\Eps^+_j)$'s are irreducible and
mutually inequivalent for $j=0,2,4,\dots,12$,
due to the fusion rules $N_{j,j'}^{28}=0$
for $j,j'=0,2,4,\dots,12$.
But we obtain
\[ \lan \eta(\a^+_{14},\Eps^+_{14}),
\eta(\a^+_{14},\Eps^+_{14}) \ran= 2 \]
since $N_{14,14}^{28}=1$
Since $[\a_{14}^+]=2[\a_4^+]$, we conclude by the same
argument as used in Example \ref{E6} that  the endomorphism
$\eta(\a^+_{14},\Eps^+_{14})$ decomposes into two
mutually inequivalent irreducible endomorphisms with
equal statistical dimensions.
We conclude that the system
\be
\begin{array}{l}
\{\eta(\a^+_j,\Eps^+_j)\mid j=0,2,4,\dots 14\} \cup
\{\eta(\a^+_j,\Eps^+_j)\eta(\tau_2,\Eps^-_2)
\mid j=0,2,4,\dots 14\}\\[.4em]
\cup \,\, \{\eta(\a^+_{14},\Eps^+_{14})_1, 
\eta(\a^+_{14},\Eps^+_{14})_2,
(\eta(\a^+_{14},\Eps^+_{14})\eta(\tau_2,\Eps^-_2))_1,
(\eta(\a^+_{14},\Eps^+_{14})\eta(\tau_2,\Eps^-_2))_2\}
\end{array}
\ee
gives the entire $\Dsys$, where
$\eta(\a^+_{14},\Eps^+_{14})_1$ and
$\eta(\a^+_{14},\Eps^+_{14})_2$ are
irreducible endomorphisms arising from decomposition of
$\eta(\a^+_{14},\Eps^+_{14})$, and
$(\eta(\a^+_{14},\Eps^+_{14})\eta(\tau_2,\Eps^-_2))_1$ and
$(\eta(\a^+_{14},\Eps^+_{14})\eta(\tau_2,\Eps^-_2))_2$
are irreducible endomorphisms arising from decomposition of
$\eta(\a^+_{14},\Eps^+_{14})\eta(\tau_2,\Eps^-_2)$.
We can then draw the dual principal graph of the
subfactor $M\otimes M^\op\subset R$ as in Fig.\ \ref{dE8},
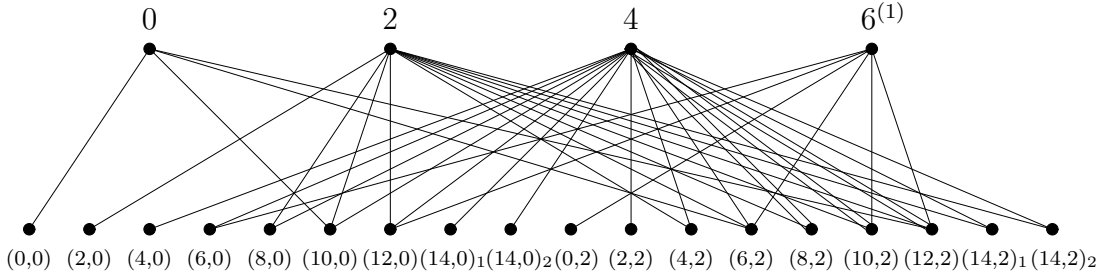
\begin{figure}[htb]
\begin{center}
\unitlength 0.8mm
\begin{picture}(190,50)
\thinlines 
\multiput(10,10)(10,0){18}{\circle*{2}}
\multiput(30,40)(40,0){4}{\circle*{2}}
\path(30,40)(10,10)
\path(30,40)(60,10)
\path(30,40)(130,10)
\path(30,40)(160,10)
\path(70,40)(20,10)
\path(70,40)(50,10)
\path(70,40)(60,10)
\path(70,40)(70,10)
\path(70,40)(120,10)
\path(70,40)(130,10)
\path(70,40)(140,10)
\path(70,40)(150,10)
\path(70,40)(160,10)
\path(70,40)(170,10)
\path(70,40)(180,10)
\path(110,40)(30,10)
\path(110,40)(40,10)
\path(110,40)(50,10)
\path(110,40)(60,10)
\path(110,40)(70,10)
\path(110,40)(80,10)
\path(110,40)(90,10)
\path(110,40)(110,10)
\path(110,40)(120,10)
\path(110,40)(130,10)
\path(109,40)(139,10)
\path(111,40)(141,10)
\path(109,40)(149,10)
\path(111,40)(151,10)
\path(109,40)(159,10)
\path(111,40)(161,10)
\path(110,40)(170,10)
\path(110,40)(180,10)
\path(150,40)(40,10)
\path(150,40)(70,10)
\path(150,40)(100,10)
\path(150,40)(130,10)
\path(150,40)(150,10)
\path(150,40)(160,10)
\put(10,5){\makebox(0,0){$\scriptstyle(0,0)$}}
\put(20,5){\makebox(0,0){$\scriptstyle(2,0)$}}
\put(30,5){\makebox(0,0){$\scriptstyle(4,0)$}}
\put(40,5){\makebox(0,0){$\scriptstyle(6,0)$}}
\put(50,5){\makebox(0,0){$\scriptstyle(8,0)$}}
\put(60,5){\makebox(0,0){$\scriptstyle(10,0)$}}
\put(70,5){\makebox(0,0){$\scriptstyle(12,0)$}}
\put(80.5,5){\makebox(0,0){$\scriptstyle(14,0)_1$}}
\put(91.5,5){\makebox(0,0){$\scriptstyle(14,0)_2$}}
\put(101,5){\makebox(0,0){$\scriptstyle(0,2)$}}
\put(110,5){\makebox(0,0){$\scriptstyle(2,2)$}}
\put(120,5){\makebox(0,0){$\scriptstyle(4,2)$}}
\put(130,5){\makebox(0,0){$\scriptstyle(6,2)$}}
\put(140,5){\makebox(0,0){$\scriptstyle(8,2)$}}
\put(150,5){\makebox(0,0){$\scriptstyle(10,2)$}}
\put(160,5){\makebox(0,0){$\scriptstyle(12,2)$}}
\put(170.5,5){\makebox(0,0){$\scriptstyle(14,2)_1$}}
\put(182,5){\makebox(0,0){$\scriptstyle(14,2)_2$}}
\put(30,45){\makebox(0,0){$0$}}
\put(70,45){\makebox(0,0){$2$}}
\put(110,45){\makebox(0,0){$4$}}
\put(152,45.5){\makebox(0,0){$6^{(1)}$}}
\end{picture}
\end{center}
\caption{The dual principal graph for the
Longo-Rehren subfactor arising from E$_8$}
\label{dE8}
\end{figure}
where we use a similar convention for labeling
vertices to the one in Fig.\ \ref{dE6}.
The procedure to get the labels for the $R$-$R$ morphisms
here is again an orbifold procedure of order 2 for the labels
$(j,\ell)$ with $=0,2,4,\dots,28$, $\ell=0,2$ with symmetry
$(j,\ell)\leftrightarrow (28-j,\ell)$.}
\end{example}

\begin{example}\labl{E(8)}{\rm 
We next study the subfactor $N\subset M$ arising from
the conformal inclusion $\SUd_5\subset\mathit{SU}(6)_1$,
as treated in \cite[Sect.\ 2.3 (iv)]{BE2}.
The 21 irreducible endomorphisms in the full $\SUd_5$
system are labelled as
$\NXN=\{\la_{(p,q)} \mid 0\le q \le p \le 5 \}$ as usual. 
Those in the ambichiral system $\MXMo$ are
labelled with the six circled vertices of the graph
$\cE^{(8)}$ in \cite[Fig.\ 11]{BE2}
and obey $\bbZ_6$ fusion rules. We label them as $\tau_\ell$,
$\ell=0,1,\dots,5$ such that the fusion rules read
$[\tau_\ell][\tau_{\ell'}]
=[\tau_{\ell+\ell' ({\mathrm{mod}}\ 6)}]$.
The non-vanishing branching coefficients
$b^+_{\tau_\ell,\la_{(p,q)}}=b_{\ell,(p,q)}$ are
\[ \begin{array}{l}
b_{0,(0,0)}=b_{0,(4,2)}=b_{1,(2,0)}=b_{1,(5,3)}
=b_{2,(3,1)}=b_{2,(5,5)}=\\[.4em]
=b_{3,(3,0)}=b_{3,(3,3)}=b_{4,(3,2)}=b_{4,(5,0)}
=b_{5,(2,2)}=b_{5,(5,2)}=1 \,
\end{array}\]
The colour zero subsystem is given by
\[ \NYN=\{\la_{(0,0)}, 
\la_{(3,0)}, \la_{(2,1)}, \la_{(5,1)}, \la_{(4,2)},
\la_{(3,3)}, \la_{(5,4)} \} \,. \]
The situation is particularly simple as this
system is still non-degenerate, i.e.\ $\NYNd=\{\id\}$.
(Note that a degenerate subsystem must be the dual of a group.)
Then $\MYMp$ consists of four endomorphisms
$\a^+_{(0,0)}$, 
$\a^+_{(5,1)}$, 
$\a^+_{(5,4)}$, 
$\a^{+,(1)}_{(3,0)}$ labelled as in \cite[Fig.\ 11]{BE2},
where we write $\a^+_{(p,q)}$ for $\a^+_{\la_{(p,q)}}$.
We study the Longo-Rehren subfactor
$M\otimes M^\op\subset R(\sys)$ arising from this
system $\sys=\MYMp$.
Again, Proposition \ref{upsglob}
implies that we only need to consider $\eta$-extensions of
color zero morphisms to obtain the entire $\Dsys$.
But in fact, due to the non-degeneracy Corollary \ref{double2}
applies and yields that $\Dsys$ is equivalent
(as $C^*$-tensor categories) to
$\NYN\times\bbZ_2$ as the subsystem of $\MYMo\subset\MXMo$
of color zero ambichirals consists of $\tau_0, \tau_3$,
obeying $\bbZ_2$ fusion rules. In particular,
here is no orbifold procedure. Alternatively,
one checks by Theorem \ref{mainirc} easily that
$\{\eta(\a_\la^+,\Eps^+_\la)\eta(\tau_k,\Eps^-_k)\mid
\la\in\NYN\,,\, k=0,3\}$ constitutes as set of 14
irreducible, mutually inequivalent endomorphisms,
hence yielding the entire quantum double system $\Dsys$.
The subfactor $\a^+_{(1,0)}(M)\subset M$
is a natural analogue of the subfactors with principal graphs
E$_6$ or E$_8$, and our Longo-Rehren subfactor corresponds to
the asymptotic inclusion of (the corresponding
hyperfinite II$_1$ subfactor of) this
inclusion. From \cite[Fig.\ 11]{BE2},
it is easy to extract the dual principal graph of the subfactor
$\a^+_{(1,0)}(M)\subset M$, drawn in Fig.\ \ref{pE(8)}.
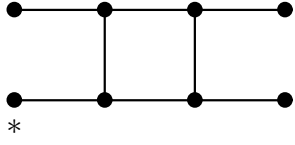
\begin{figure}[htb]
\begin{center}
\unitlength 1.2mm
\begin{picture}(50,20)
\thicklines 
\multiput(10,5)(10,0){4}{\circle*{1.5}}
\multiput(10,15)(10,0){4}{\circle*{1.5}}
\path(10,5)(40,5)
\path(10,15)(40,15)
\path(20,5)(20,15)
\path(30,5)(30,15)
\put(10,2){\makebox(0,0){$*$}}
\end{picture}
\end{center}
\caption{The (dual) principal graph for the subfactor
$\a^+_{(1,0)}(M)\subset M$}
\label{pE(8)}
\end{figure}
Using Proposition \ref{dualpg} it is now a straight-forward
calculation yielding the dual principal graph of the
associated Longo-Rehren inclusion, displayed in
Fig.\ \ref{dcE(8)}.
\begin{figure}[htb]
\begin{center}
\unitlength 1.0mm
\begin{picture}(150,40)
\thinlines 
\multiput(10,10)(10,0){14}{\circle*{2}}
\multiput(30,30)(30,0){4}{\circle*{2}}
\path(30,30)(10,10)
\path(30,30)(20,10)
\path(30,30)(30,10)
\path(30,30)(60,10)
\path(60,30)(40,10)
\path(60,30)(50,10)
\path(60,30)(60,10)
\path(60,30)(70,10)
\path(60,30)(80,10)
\path(60,30)(90,10)
\path(90,30)(60,10)
\path(90,30)(70,10)
\path(90,30)(80,10)
\path(90,30)(90,10)
\path(90,30)(100,10)
\path(90,30)(110,10)
\path(120,30)(90,10)
\path(120,30)(120,10)
\path(120,30)(130,10)
\path(120,30)(140,10)
\put(10,5){\makebox(0,0){$\scriptstyle(0,0;0)$}}
\put(20,5){\makebox(0,0){$\scriptstyle(3,0;3)$}}
\put(30,5){\makebox(0,0){$\scriptstyle(3,3;3)$}}
\put(40,5){\makebox(0,0){$\scriptstyle(5,1;0)$}}
\put(50,5){\makebox(0,0){$\scriptstyle(5,4;3)$}}
\put(60,5){\makebox(0,0){$\scriptstyle(4,2;0)$}}
\put(70,5){\makebox(0,0){$\scriptstyle(2,1;0)$}}
\put(80,5){\makebox(0,0){$\scriptstyle(2,1;3)$}}
\put(90,5){\makebox(0,0){$\scriptstyle(4,2;3)$}}
\put(100,5){\makebox(0,0){$\scriptstyle(5,4;0)$}}
\put(110,5){\makebox(0,0){$\scriptstyle(5,1;3)$}}
\put(120,5){\makebox(0,0){$\scriptstyle(3,0;0)$}}
\put(130,5){\makebox(0,0){$\scriptstyle(3,3;3)$}}
\put(140,5){\makebox(0,0){$\scriptstyle(0,0;3)$}}
\put(30,35){\makebox(0,0){$(0,0)$}}
\put(60,35){\makebox(0,0){$(5,1)$}}
\put(90,35){\makebox(0,0){$(5,4)$}}
\put(122,35.3){\makebox(0,0){$(3,0)^{(1)}$}}
\end{picture}
\end{center}
\caption{The dual principal graph for the
Longo-Rehren subfactor arising from $\cE^{(8)}$}
\label{dcE(8)}
\end{figure}
Here we used the short-hand notation $(p,q;\ell)$ for
$\eta(\a^+_{(p,q)},\Eps^+_{(p,q)})\eta(\tau_\ell,\Eps^-_\ell)$
and an obvious notation for the vertices labelled
by morphisms in $\MYMp$.
It seems that the system $\Dsys$ is equivalent to the system
of $Q_\infty$-$Q_\infty$ bimodules arising from the
asymptotic inclusion $Q\vee (Q'\cap Q_\infty)$ of the
hyperfinite II$_1$ subfactor $P\subset Q$ with principal
graph A$_7$ of Jones \cite{J}, but we have no proof.}
\end{example}

\begin{example}\labl{E(12)}{\rm 
We next study the subfactor $N\subset M$ arising from
the conformal inclusion $SU(3)_9\subset ($E$_6)_1$,
as treated in \cite[Sect.\ 6.4]{BE3}.
The 55 irreducible endomorphisms in $\NXN$ are labelled with
$\la_{(p,q)}$, $0\le q \le p \le 9$ as usual.
The chiral system $\MXMp$ corresponds to the vertices
of the graph ${\mathcal E}^{(12)}$, and the ambichiral
system $\MXMo$ to the three vertices marked with circles
in \cite[Fig.\ 12]{BE3}, obeying the $\bbZ_3$ fusion rules.
We label them as $\tau_\ell$, $\ell=0,1,2$, so that
$[\tau_\ell][\tau_{\ell'}]=[\tau_{\ell+\ell'({\mathrm{mod}}\ 3)}]$.

The non-vanishing branching coefficients
$b^+_{\tau_\ell,\la_{(p,q)}}=b_{\ell,(p,q)}$ are
\[ \begin{array}{l}
b_{0,(0,0)}=b_{0,(5,1)}=b_{0,(5,4)}=b_{0,(8,4)}
=b_{0,(9,0)}=b_{0,(9,9)}=\\[.4em]
=b_{1,(4,2)}=b_{1,(7,1)}=b_{1,(7,7)}
=b_{2,(4,2)}=b_{2,(7,1)}=b_{2,(7,7)}=1 \,
\end{array}\]

The color zero subsystem $\NYN$ is given by those 19 morphisms
$\la_{(p,q)}\in\NXN$ subject to $p+q\in 3\bbZ$.
It now contains the simple currents, and as a
consequence \cite[Lemma 6.11]{BE3}
we have $\NYNd=\{\la_{(0,0)},\la_{(9,0)},\la_{(9,9)}\}$.
Then the system $\MYMp$ consists of four endomorphisms
$\a^+_{(0,0)}=\tau_0$, $\a^+_{(2,1)}$, $\tau_1$, $\tau_2$,
so that in particular $\MXMo=\MYMo$.
(We use labels as in \cite[Fig.\ 12]{BE3}, apart from
denoting the ambichiral $\eta_j$ by $\tau_j$ here, $j=1,2$,
as $\eta_j$ is obviously no suitable notation when
considering $\eta$-extensions.)
As usual, we study the Longo-Rehren subfactor
$M\otimes M^\op\subset R(\sys)$ arising from $\sys=\MYMp$.
Since the degenerate morphisms appear in the vacuum column of the
modular invariant, we only need to consider $\eta$-extensions of
$\a^+_\la$ and $\tau$ with $\la\in\NYN$ and $\tau\in\MYMo$ only,
thanks to Proposition \ref{upsglob}.
The subfactor $\a^+_{(1,0)}(M)\subset M$ is again a natural
analogue of the subfactors with principal graphs E$_6$ and
E$_8$. From \cite[Fig.\ 12]{BE3} it is easy extract the dual
principal graph of the subfactor $\a^+_{(1,0)}(M)\subset M$,
drawn in Fig.\ \ref{pE(12)}.
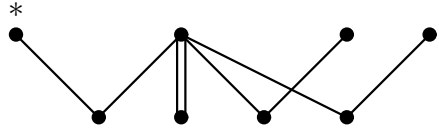
\begin{figure}[htb]
\begin{center}
\unitlength 1.1mm
\begin{picture}(70,20)
\thicklines 
\multiput(20,5)(10,0){4}{\circle*{1.5}}
\multiput(10,15)(20,0){3}{\circle*{1.5}}
\put(60,15){\circle*{1.5}}
\path(10,15)(20,5)
\path(20,5)(30,15)
\path(30.5,5)(30.5,15)
\path(29.5,5)(29.5,15)
\path(30,15)(40,5)
\path(50,15)(40,5)
\path(30,15)(50,5)
\path(60,15)(50,5)
\put(10,18){\makebox(0,0){$*$}}
\end{picture}
\end{center}
\caption{The (dual) principal graph for the
subfactor $\a^+_{(1,0)}(M)\subset M$}
\label{pE(12)}
\end{figure}
By Lemma \ref{etalpmpc} we find that
$\eta(\a^+_{(0,0)},\Eps^+_{(0,0)})$,
$\eta(\a^+_{(2,1)},\Eps^+_{(2,1)})$,
$\eta(\a^+_{(3,0)},\Eps^+_{(3,0)})$,
$\eta(\a^+_{(3,3)},\Eps^+_{(3,3)})$,
$\eta(\a^+_{(4,2)},\Eps^+_{(4,2)})$, and
$\eta(\a^+_{(5,1)},\Eps^+_{(5,1)})$ are irreducible
and mutually inequivalent endomorphisms of $R$.
We similarly obtain
\[ \lan \eta(\a^+_{(6,3)},\Eps^+_{(6,3)}),
\eta(\a^+_{(6,3)},\Eps^+_{(6,3)}) \ran =3 \,, \]
and find that $\eta(\a^+_{(6,3)},\Eps^+_{(6,3)})$ is
disjoint from the others since $\la_{(6,3)}$ is a fixed point
of the simple currents. Consequently, the decomposition of
$[\eta(\a^+_{(6,3)},\Eps^+_{(6,3)})]$ yields three new
irreducible sectors.
Next, Theorem \ref{mainirc} yields
$\lan\eta(\a^+_\la,\Eps^+_\la),\eta(\tau_\ell,\Eps^-_\ell)\ran=0$
for $\la\in\NYN$ and $\ell=1,2$, because $[\tau_1]$ and
$[\tau_2]$ do not appear as subsectors of $[\a^+_\la]$
for $\la\in\NYNd$.
Similarly we find that
$\eta(\a^+_\la,\Eps^+_\la)\eta(\tau_\ell,\Eps^-_\ell)$
are disjoint for different $\ell=0,1,2$.
We now have a set of irreducible, mutually inequivalent
endomorphisms of $R$ consisting of 27 endomorphisms,
which can be considered as $\NYN/\bbZ_3 \times \bbZ_3$.
Here for the orbifold $\NYN/\bbZ_3$, 18 objects collapse into
6 objects by identification arising from a $\bbZ_3$ symmetry,
and the fixed point of the symmetry splits into 3 objects.
The total number of the irreducible objects is therefore
$(6+3)\times 3=27$.
Let the statistical dimensions of the irreducible morphisms
appearing in the decomposition of 
$\eta(\a^+_{(6,3)},\Eps^+_{(6,3)})$ be $d_1,d_2, d_3$ respectively.
We then have $d_1+d_2+d_3=d_{(6,3)}$.  The square sum
$d_1^2+d_2^2+d_3^2$ attains the minimum $d_{(6,3)}^2/3$ with
$d_1=d_2=d_3=d_{(6,3)}/3$ under the constraint
$d_1+d_2+d_3=d_{(6,3)}$. Assume for contradiction
that we are off the minimum.
Then $d_1^2+d_2^2+d_3^2 > d_{(6,3)}^2/3$, and in turn
the global index of the system  $\NYN/\bbZ_3 \times \bbZ_3$
is strictly bigger than $[[\NYN]]$.
But then it exceeds $[[\Dsys]]=[[\MYMp]]^2$
because \cite[Prop.\ 3.1]{BE4} tells us that
$[[\MYMp]]^2=[[\NYN]][[\MYMo]]/\sum_{\la\in\NYNd}d_\la Z_{\la,0}$
which obviously yields $[[\MYMp]]^2=[[\NYN]]$ here;
contradiction. We conclude that
$d_1=d_2=d_3=d_{(6,3)}/3$ and that the above set
of morphisms gives the entire system $\Dsys$.}
\end{example}

\begin{example}\labl{E(24)}{\rm 
We finally study the subfactor $N\subset M$ arising from
the conformal inclusion $SU(3)_{21}\subset ($E$_7)_1$,
as treated in \cite[Sect.\ 6.4]{BE3}.
The 253 irreducible endomorphisms in $\NXN$ are labelled with
$\la_{(p,q)}$, $0\le q \le p \le 21$, as usual.
The morphisms in $\MXMp$ correspond to the vertices
of the graph ${\mathcal E}^{(24)}$, and those in the
ambichiral system $\MXMo$ to the two encircled vertices
in \cite[Fig.\ 13]{BE3}, the latter
obeying the $\bbZ_2$ fusion rules. 
We label them as $\tau_\ell$ here, $\ell=0,1$,
so that $[\tau_1][\tau_1]=[\tau_0]$.
(Note that our $\tau_1$ is denoted
by $\epsilon$ in\cite[Fig.\ 13]{BE3}.)
The color zero subsystem $\NYN$ is given by those 85 morphisms
$\la_{(p,q)}\in\NXN$ subject to $p+q\in 3\bbZ$.
It now contains the simple currents, and as a consequence
we have $\NYNd=\{\la_{(0,0)},\la_{(21,0)},\la_{(21,21)}\}$.
The system $\MYMp$ consists of eight endomorphisms
$\a^+_{(0,0)}=\tau_0$, $\a^+_{(2,1)}$, $\a^{+,(1)}_{(4,2)}$, 
$\a^{+,(2)}_{(4,2)}$, $\a^+_{(3,0)}$, $\a^+_{(3,3)}$, 
$\a^{+,(1)}_{(5,1)}$, and $\tau_1$.
We study the Longo-Rehren subfactor
$M\otimes M^\op\subset R(\sys)$ arising from this
system $\sys=\MYMp$, which corresponds to the asymptotic
inclusion of the subfactor $\a^+_{(1,0)}(M)\subset M$ as
the natural analogue of the subfactors with principal graphs
E$_6$, E$_8$. From \cite[Fig.\ 13]{BE3} we extract
the dual principal graph of the subfactor
$\a^+_{(1,0)}(M)\subset M$, displayed in Fig.\ \ref{pE(24)}.
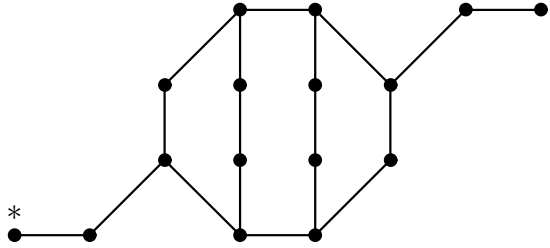
\begin{figure}[htb]
\begin{center}
\unitlength 1mm
\begin{picture}(90,40)
\thicklines 
\multiput(10,5)(10,0){2}{\circle*{1.5}}
\multiput(40,5)(10,0){2}{\circle*{1.5}}
\multiput(40,35)(10,0){2}{\circle*{1.5}}
\multiput(70,35)(10,0){2}{\circle*{1.5}}
\multiput(30,15)(10,0){4}{\circle*{1.5}}
\multiput(30,25)(10,0){4}{\circle*{1.5}}
\path(10,5)(20,5)
\path(20,5)(30,15)
\path(30,15)(30,25)
\path(30,25)(40,35)
\path(40,35)(50,35)
\path(50,35)(60,25)
\path(60,25)(60,15)
\path(60,15)(50,5)
\path(50,5)(40,5)
\path(40,5)(30,15)
\path(40,5)(40,35)
\path(50,5)(50,35)
\path(60,25)(70,35)
\path(70,35)(80,35)
\put(10,8){\makebox(0,0){$*$}}
\end{picture}
\end{center}
\caption{The (dual) principal graph for the
subfactor $\a^+_{(1,0)}(M)\subset M$}
\label{pE(24)}
\end{figure}
Since the degenerate morphisms appear in the vacuum
column of the modular invariant, the situation is similar
to Examples \ref{E8} and \ref{E(12)}.
Their $\bbZ_3$ symmetry has $\la_{(14,7)}\in\NYN$ as a fixed point,
and the other 84 endomorphisms give 28 orbits under this symmetry.
We then have
$\lan \eta(\a^+_{(14,7)},\Eps^+_{(14,7)}),
\eta(\a^+_{(14,7)},\Eps^+_{(14,7)})\ran=3$.
Along the same lines as in Example \ref{E(12)},
we conclude that the system $\Dsys$ contains
$(28+3)\times 2=62$ irreducible endomorphisms,
corresponding to $\NYN/\bbZ_3\times\bbZ_2$.
Namely, the 28 irreducible endomorphisms
$\eta(\a^+_{\la},\Eps^+_\la)$ where we select one
$\la\in\NYN$ of each $\bbZ_3$ orbit together with
the three irreducible endomorphisms of equal
statistical dimensions arising from decomposition of
$\eta(\a^+_{(14,7)},\Eps^+_{(14,7)})$
correspond to $\NYN/\bbZ_3$, and the blowing up by
$\bbZ_2$ arises from multiplication with
$\eta(\tau_\ell,\Eps^-_\ell)$, $\ell=0,1$.}
\end{example}


\vspace{0.5cm}
\begin{footnotesize}
\noindent{\it Acknowledgment.}
Part of this work was done during visits of the third
author to the University of Wales Cardiff.
We thank M. Izumi and K.-H. Rehren for helpful discussions on
\cite{I1,I2}, \cite{R7}, and we are grateful to J.E. Roberts
and S. Yamagami for helpful comments on category theory.
Y.K. thanks A. Ocneanu for
discussions at the Taniguchi Conference at Nara, Japan,
in December 1998 concerning his lectures given there.
We gratefully acknowledge the financial support of
EPSRC (U.K.), the EU TMR Network in Non-Commutative Geometry,
Grant-in-Aid for Scientific Research, Ministry of
Education (Japan), the Mitsubishi Foundation,
University of Tokyo, and the University of Wales.
\end{footnotesize}

\begin{footnotesize}

\end{footnotesize}
\end{document}